%% file: example.tex
\documentclass[submission,copyright,creativecommons]{eptcs}
\pdfoutput=1

\usepackage{paralist}
\usepackage{amsmath}
\usepackage{underscore}           %
\usepackage{xstring}
\usepackage{amsthm}
\usepackage{floatrow}
\usepackage{hyperref}
\usepackage{amsmath}
\usepackage{amsthm}
\usepackage{amssymb} %
\usepackage{graphicx}
\usepackage{float}
\def\sagdir{sag}
\def\rootdir{.}
\input{utils/preamble_act}
\makeatletter
\def\mathcolor#1#{\@mathcolor{#1}}
\def\@mathcolor#1#2#3{%
  \protect\leavevmode
  \begingroup
    \color#1{#2}#3%
  \endgroup
}
\makeatother
\input{symbols.tex}
\newcommand{\proceedings}[2]{\ifthenelse{\boolean{proceedings}}{#1}{#2}}

\definecolor{transmuted}{RGB}{0,0,0}

\providecommand{\event}{Preprint} %

\theoremstyle{definition}

\newlength{\mysep}
\newtheorem{theorem}{Theorem}
\newtheorem*{theorem*}{Theorem}

\newtheorem*{proposition*}{Proposition}

\newtheorem*{lemma*}{Lemma}

\newtheorem{definition}[theorem]{Definition}
\newtheorem*{definition*}{Definition}

\AtBeginDocument{
        \addtolength{\abovedisplayskip}{-0.4ex}
        \addtolength{\abovedisplayshortskip}{-0.2ex}
        \addtolength{\belowdisplayskip}{-0.4ex}
        \addtolength{\belowdisplayshortskip}{-0.2ex}
}

\def\scalemed{0.55}

\title{Diagrammatic Negative Information}
    \author{Vincent Abbott
    \institute{College of Engineering, Computing and Cybernetics\\
    Australian National University\\
    Canberra, ACT, Australia}
    \email{Vincent.Abbott@anu.edu.au}
    \and
    Gioele Zardini
    \institute{
    Laboratory for Information and Decision Systems\\
    Massachusetts Institute of Technology\\
    Cambridge, MA, USA}
    \email{gzardini@mit.edu}
}
\def\titlerunning{Diagrammatic Negative Information}
\def\authorrunning{V. Abbott, G. Zardini}

\begin{document}
\maketitle

\begin{abstract}
The flow of information through a complex system can be readily understood with category theory. 
However, negative information (e.g., what is \textit{not }possible) does not have an immediately evident categorical representation. 
The formalization of nategories using unconventional composition addresses this issue, and lets imposed limitations on categories be considered. 
However, traditional nategories abandon core categorical constructs and rely on extensive mathematical development. This creates a divide between the consideration of positive and negative information composition. 
In this work, we show that negative information can be considered in a natural categorical manner. This is aided by functor string diagrams, a novel flexible diagrammatic approach that can intuitively show the operation of hom-functors and natural transformations in expressions. This insight reveals how to consider the composition of negative information with foundational categorical constructs without relying on enrichment.
We present diagrammatic means to consider not only nategories, but preorders more broadly. 
This paper introduces diagrammatic methods for the consideration of triangle inequalities and n  co-designs $\displaystyle \mathbf{DP/Feas_{Bool}}$, showing how important cases of negative information composition can be categorically and diagrammatically approached. 
In particular, we develop systematic tools to rigorously consider imposed limitations on systems, advancing our mathematical understanding, and present intuitive diagrams which motivate widespread adoption and usage for various applications.
\end{abstract}

\section{Introduction}
Category theory is the formal study of composition. 
Morphisms express how components can be connected, and thus represent various ``paths'' which can be taken. 
Morphisms show both that objects can be linked and, given the multitude morphisms possible, associate these links with information.

However, practical applications often need to consider negative (complementary) information (e.g., paths that are not possible).
In trying to design or plan a connection between certain states, it is important to consider minimality, and impossibilities from which to work. This concept of negative information is central to widely used techniques such as Dijkstra's algorithm or A* search \cite{delling-engineering-2009}.

However, negative information has a peculiar composition structure. 
A limitation on~$\displaystyle a\rightarrow b$ carries through to paths~$\displaystyle b\rightarrow c$ and $\displaystyle a\rightarrow c$, which does not match typical categorical structure. 
This follows from paths $\displaystyle b\rightarrow c$ and $\displaystyle a\rightarrow c$ being composable into $\displaystyle a\rightarrow b$ routes, hence, a ban on the former imposes itself on the latter. 
This is not in accordance with standard category theory which has $\displaystyle a\rightarrow b$ compose with morphisms to $\displaystyle a$ or from $\displaystyle b$.

Nategories, recently introduced in~\cite{zarACT2022}, offer a means to consider these effects by introducing additional mathematical structure to categories, abandoning the general rule that morphism + morphism $\rightarrow$ morphism and adding ``norphisms'', representing bans, which follow morphism + norphism $\rightarrow$ norphism, representing carry through effects. 
Certain norphisms can be considered purely categorical using an enriched construct on the mathematically advanced $\displaystyle \mathbf{\textcolor[rgb]{0.29,0.56,0.89}{P}\textcolor[rgb]{0.82,0.01,0.11}{N}}$ (positive-negative) category.

Diagrammatic category theory maps categorical properties onto diagrams, associating algebraic rules of categories with the manipulation of diagrams \cite{selinger-survey-2011}. 
This is central to monoidal string diagrams, whose categorical structure corresponds exactly to diagrams with graphical isomorphism rules. However, traditional monoidal string diagrams are limited to expressing simple connections between objects and morphisms, and are not able to readily represent functors and natural transformations, which requires further constructs \cite{marsden-category-2014, nakahira-diagrammatic-2023, piedeleu-introduction-2023}.

String diagrams developed by Marsden and Nakahira \cite{marsden-category-2014,nakahira-diagrammatic-2023,piedeleu-introduction-2023} represent categories with colored regions, associating functors to the boundaries of regions and natural transformations to the intersection of boundaries. This graphically encodes the behavior of functors and natural transformations through similar intuitive isotopy rules. However, they are designed with to focus on algebraic constructs such as monoids, Kan extensions, and adjunctions, but lack the flexibility to be used for or to give insight for applied category theory.

Functor string diagrams are a novel approach which streamline string diagrams using solid principles which ensure that complex diagrams are decipherable \cite{FSD-preprint, abbott-robust-2023}. They are used as the mathematical foundations of neural circuit diagrams for machine learning \cite{abbott-neural-2023}. However, they can also be used to consider abstract algebraic ideas in an intuitive manner and, as we will soon show, the manipulation of natural transformations and hom-functors they encourage are the key to revealing how foundational category theory can encompass norphisms and negative information.
\section{Background}
\subsection{Nategories}
The current means of considering negative information in a category is with \textit{nategories}. Nategories are categories embedded with additional elements, called norphisms, with atypical composition rules. 
\begin{definition}[Nategory]\label{def:nategory}
    A locally small \emph{nategory} \CatC is a locally small category with the following additional structure.
    For each pair of objects $\Obja,\Objb \setin \Ob_{\CatC}$, in addition to the set of morphisms $\HomSet{\CatC}{\Obja}{\Objb}$, we also specify:
    \begin{compactitem}
        \item A set of norphisms $\NomSet{\CatC}{\Obja}{\Objb}$.
        \item An \emph{incompatibility relation}, which we write as a binary function
              \begin{equation} \label{eq:010-nategory-incompat}
                  \nmincompat{\Obja}{\Objb}\colon\NomSet{\CatC}{\Obja}{\Objb} \cartprod \HomSet{\CatC}{\Obja}{\Objb}   \to 2.
              \end{equation}
    \end{compactitem}
    For all triples $\Obja,\Objb,\Objc$, in addition
    to the morphism composition function
    \begin{equation} \label{eq:011-nategory-comp}
        \mthen_{\Obja\Objb\Objc} \colon \HomSet{\CatC}{\Obja}{\Objb} \cartprod \HomSet{\CatC}{\Objb}{\Objc} \to \HomSet{\CatC}{\Obja}{\Objc},
    \end{equation}
    we require the existence of two \textbf{inexact }norphism composition functions
    \begin{equation}\label{eq:012-nategory-norcomp}
        \begin{aligned}
            \northenbsymb_{\Obja\Objb\Objc} & \colon \HomSet{\CatC}{\Obja}{\Objb} \cartprod \NomSet{\CatC}{\Obja}{\Objc} \to \NomSet{\CatC}{\Objb}{\Objc},
            \\
            \northenasymb_{\Obja\Objb\Objc} & \colon \NomSet{\CatC}{\Obja}{\Objc} \cartprod \HomSet{\CatC}{\Objb}{\Objc} \to \NomSet{\CatC}{\Obja}{\Objb},
        \end{aligned}
    \end{equation}
    and we ask that they satisfy two ``equivariance'' conditions:
    \begin{align}
        \nmincompat{\Objb}{\Objc}(\northenb{\mora}{\nora}, \morb)
         & \Rightarrow
        \nmincompat{\Obja}{\Objc} ( \nora, \morab )\label{eq:cond1}\tag{equiv-1}, \\
        \nmincompat{\Obja}{\Objb} (\northena{\nora}{\morb}, \mora)
         & \Rightarrow
        \nmincompat{\Obja}{\Objc} ( \nora, \morab )\label{eq:cond2}\tag{equiv-2}.
    \end{align}
\end{definition}

As can be seen, this construction extends typical categories, abandoning the universal notion of composition in order to consider negative information. However, using a novel diagrammatic approach, we can show these two forms of composition to be the natural interactions between hom-functors and functions $\displaystyle \mathbf{C}( a,b)\rightarrow 2$ in the category $\displaystyle \mathbf{Set}$.

\subsection{Functor String Diagrams}
\begin{figure}[htb]
    \floatbox[{\capbeside\thisfloatsetup{capbesideposition={left,top}}}]{figure}[\FBwidth]
    {\caption{Functor string diagrams show categorical expressions by composing vertical sections, each of which represents an object or a morphism.}
    \label{fig:FSD}}
    {\includegraphics[scale=\scalemed]{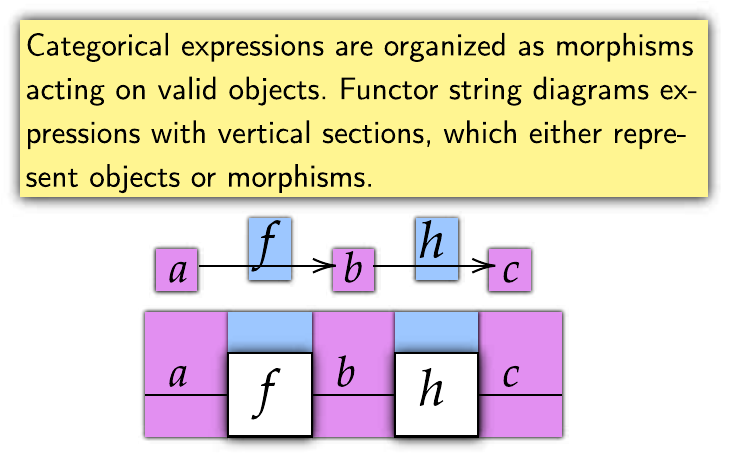}}
\end{figure}

To understand these constructs better, we use functor string diagrams, a novel diagrammatic method tailored to expressing the behavior of functors and natural transformations within a category \cite{FSD-preprint, abbott-robust-2023}. 
We show these in Figure \ref{fig:FSD}. Diagrams are partitioned into \emph{vertical} sections, each associated to an object or morphism. 
Functors are represented by functor wires passing over relevant objects and morphisms, which intuitively modify the underlying objects or morphisms while encoding the preservation of composition.
\FloatBarrier

\begin{figure}[htb] \floatbox[{\capbeside\thisfloatsetup{capbesideposition={left,top}}}]{figure}[\FBwidth]
    {\caption{For more complex expressions, we develop equivalent expressions which lets categorical algebra be understood with graphical intuition.}
    \label{fig:FSDhom}}
    {\includegraphics[scale=\scalemed]{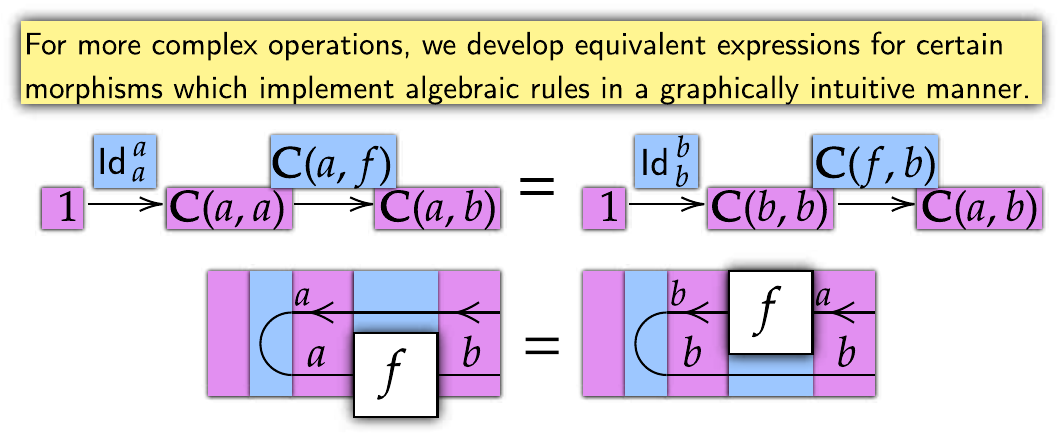}}
\end{figure}

Hom-functors are represented by a wire with a leftward arrow labeled with the relevant object. Natural transformations between hom-functors $\displaystyle \mathbf{C}( c,\_)\rightarrow \mathbf{C}( a,\_)$ correspond to morphisms $\displaystyle f:a\rightarrow c$ by the Yoneda lemma. As shown in Figure \ref{fig:FSDhom}, natural transformations are drawn on the functor wire, encoding their ability to pass over underlying morphisms while maintaining equivalence.

\FloatBarrier

Consider that norphisms are defined by their interaction with $\nmincompat{\Obja}{\Objb}\colon\NomSet{\CatC}{\Obja}{\Objb} \cartprod \HomSet{\CatC}{\Obja}{\Objb}   \to 2$. Currying $\nmincompat{\Obja}{\Objb}$, we see that a norphism $\nora \in \NomSet{\CatC}{\Obja}{\Objb}$ corresponds to a function $\nmincompat{}{} \nora:\HomSet{\CatC}{\Obja}{\Objb}\rightarrow 2$. This lets us consider the curried form of norphisms, $\nmincompat{}{} \nora:\HomSet{\CatC}{\Obja}{\Objb}\rightarrow 2$, as morphisms $\textbf{C}(a,b)\rightarrow 2$ in the category $\textbf{Set}$.

This reveals two manners that morphisms can interact with $\displaystyle \mathbf{C}( a,b)\rightarrow 2$, as morphisms $\displaystyle \mathbf{C}( c,b)$ with a $\displaystyle \mathbf{C}( a,\_)$ functor applied, or natural transformations mapping $\displaystyle \mathbf{C}( c,\_)\rightarrow \mathbf{C}( a,\_)$ over an underlying $\displaystyle b$ object wire.

\FloatBarrier
\begin{figure}[htb] \floatbox[{\capbeside\thisfloatsetup{capbesideposition={left,top}}}]{figure}[\FBwidth]
    {\caption{We can re-express norphism and their exact composition with morphisms using functor string diagrams.}
    \label{fig:FSDnorph}}
    {\includegraphics[scale=\scalemed]{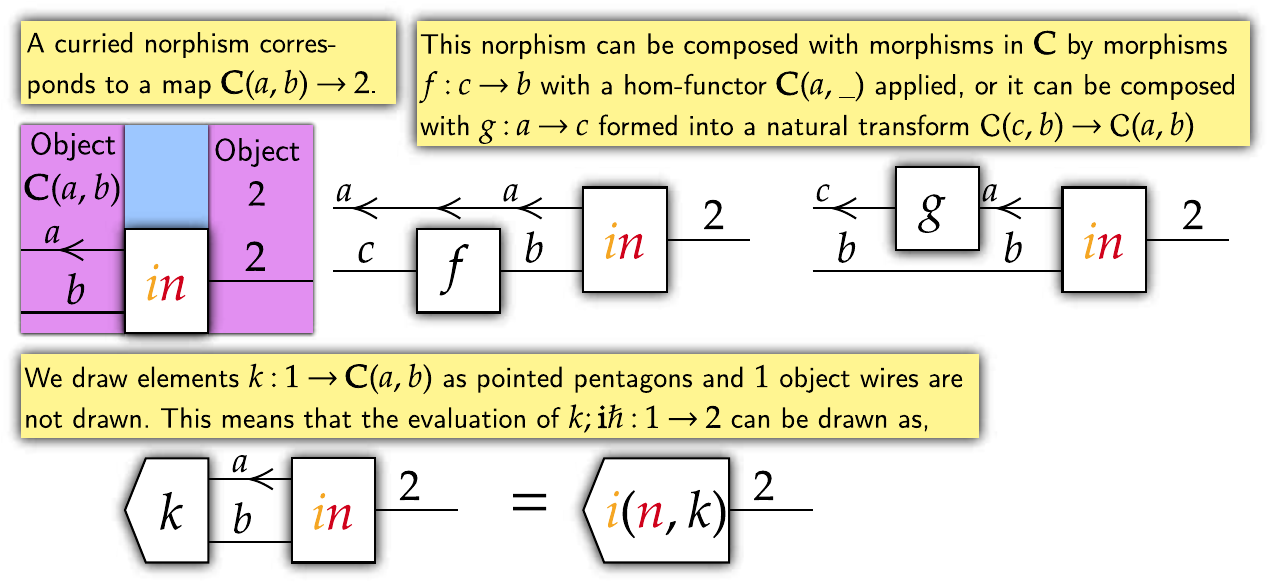}}
\end{figure}

The shape of inexact nategorical compositions are the same as the above. Inexact nategory composition yields $\northenbsymb_{\Obja\Objb\Objc} \colon \HomSet{\CatC}{\Obja}{\Objb} \cartprod \NomSet{\CatC}{\Obja}{\Objc} \to \NomSet{\CatC}{\Objb}{\Objc}$ and $\northenasymb_{\Obja\Objb\Objc} \colon \NomSet{\CatC}{\Obja}{\Objc} \cartprod \HomSet{\CatC}{\Objb}{\Objc} \to \NomSet{\CatC}{\Obja}{\Objb}$. To uniquely express these, we require reference to the constituent morphism (morphism in $\CatC$) and norphism (corresponding to a function $\HomSet{\CatC}{\Obja}{\Objb}\rightarrow 2$). We do this by copying the shape of the above exact compositions, but draw the norphism composition symbols on the connecting wires, indicating that we are not perofrming a typical composition. The equivariance condition, then, is implemented by stating that these inexact compositions imply the exact compositions.

\begin{figure}[htb] \floatbox[{\capbeside\thisfloatsetup{capbesideposition={left,top}}}]{figure}[\FBwidth]
    {\caption{We use dotted thick wires to indicate morphisms produced from inexact norphism composition. Note that inexact composition implies exact composition.}
    \label{fig:inexact-composition}}
    {\includegraphics[scale=\scalemed]{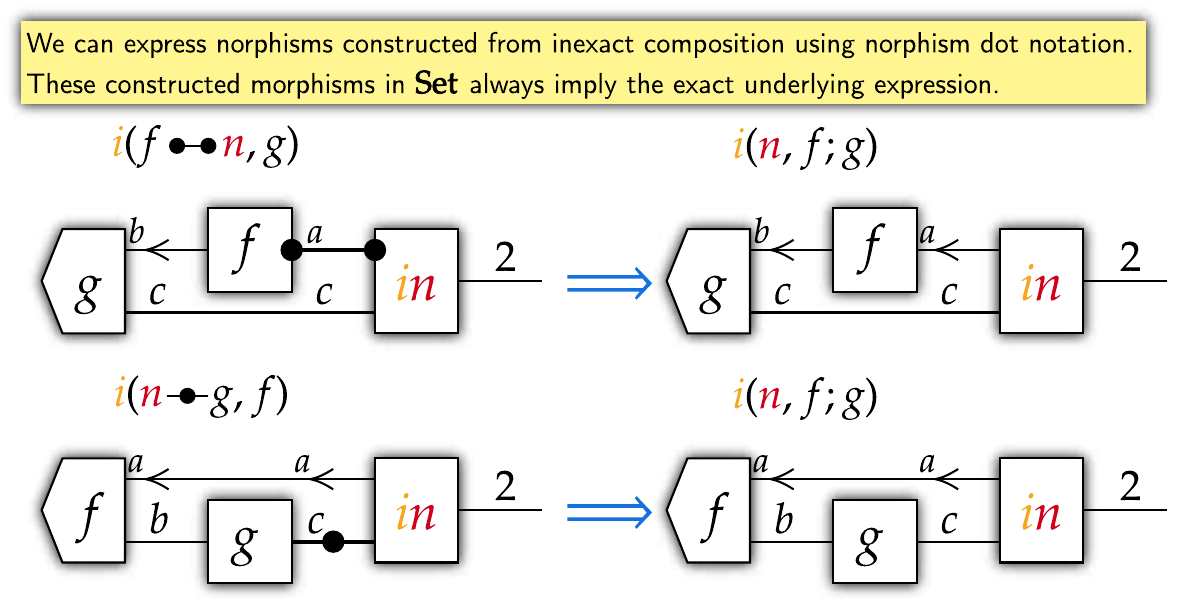}}
\end{figure}

So far, we have developed diagrammatic tools to represent nategories. However, we have relied on implications instead of pure equivalence of expressions, meaning the application of categorical algebra is unclear. In the next section, we will provide tools to diagrammatically view preorders as transforms which are pure equivalences. This perspective allows preorders to be understood as equivalences, avoiding the complexities of $2$-category algebra.

\section{Diagrammatic Preorders}
Preorders have second-order categorical structure. This complicates the mathematics, and makes expressions and algebra confusing to understand. This complexity obfuscates the insights preorders offer and hampers their broader application despite their utility for understanding design and trade-offs in engineered systems among other fields (see Section \ref{sec:diagrammatic-co-designs}).

Here, we introduce a diagrammatic means of considering preorders and $\displaystyle 2$-categories more generally with functor string diagrams, letting them be simultaneously be considered along with the hom-functors and natural transformations we will frequently use. This makes the insights they offer easier to perceive. This diagrammatic framework intuitively captures norphisms and co-designs $\displaystyle \mathbf{DP/Feas_{Bool}}$ in a straightforward manner. This both aids mathematical rigor by motivating the use of theoretical tools indicated by diagrams and encourages adoption for applied cases given the clear manner in which diagrams enable ideas to be expressed.

\paragraph{Preorders.} Preorders $\displaystyle \mathbf{X}_{\leqslant}$ are a generalization of $\displaystyle \leqslant $ relationships. A set $\displaystyle X$ with internal preorders can be expressed as a category $\displaystyle \mathbf{X}_{\leqslant }$ where objects are elements and morphisms are $\displaystyle \leqslant $ relationships. This encodes for the fact that elements exhibit $\displaystyle x\leqslant _{x}^{x} \ x$, the identity relationship, and that $\displaystyle x\leqslant _{y}^{x} \ y$ and $\displaystyle y\leqslant _{z}^{y} \ \ z$ implies there is a relationship $\displaystyle x\leqslant _{z}^{x} \ z$. These correspond to morphisms $\displaystyle \mathbf{X}_{\leqslant }( x,x)$ and composition $\displaystyle \mathbf{X}_{\leqslant }( x,y) \times \mathbf{X}_{\leqslant }( y,z)\rightarrow \mathbf{X}_{\leqslant }( x,z)$.

\paragraph{Preorders as Transforms.} To maintain consistent types, we introduce an initial object $\displaystyle 1$ such that each object $\displaystyle x$ of $\displaystyle \mathbf{X}_{\leqslant} \cup 1$ has a unique identifying morphism $\displaystyle x_{x}^{1} :1\rightarrow x$. Then, $\displaystyle \leqslant _{y}^{x}$ relationships let us compose $y_{y}^{1} = \displaystyle x_{x}^{1} ;\leqslant _{y}^{x}$, mapping from the unique identifier for $\displaystyle x$ to the unique identifier for $\displaystyle y$.
\begin{figure}[h!] \floatbox[{\capbeside\thisfloatsetup{capbesideposition={left,top}}}]{figure}[\FBwidth]
    {\caption{
A set with preorder structure can be expressed as a category where morphisms correspond to $\displaystyle \leqslant $ relationships. This allows for the transformation of objects into superior ones. With a unique identifier from an initial object $\displaystyle 1$, we can fix types and state that $\displaystyle y_{y}^{1} = x_{x}^{1} ;\leqslant _{y}^{x}$.
    }
    \label{fig:preorder}}
    {\includegraphics[scale=\scalemed]{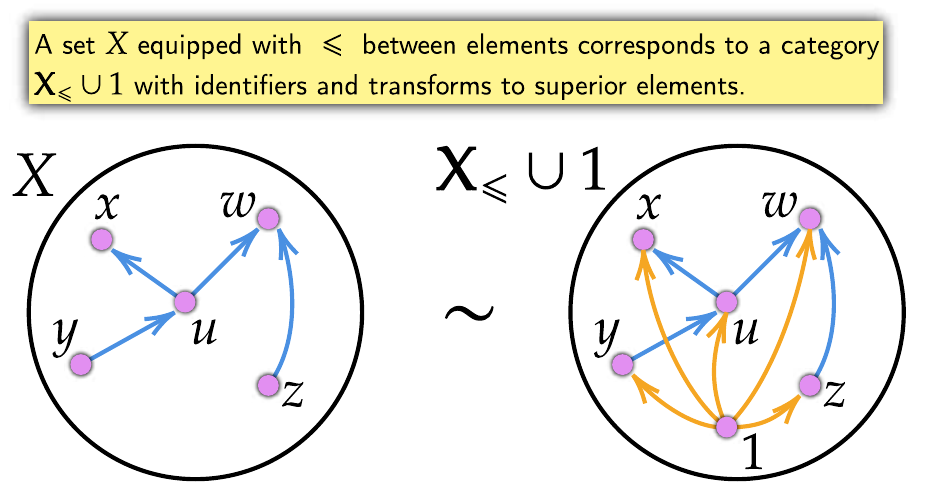}}
\end{figure}
\FloatBarrier

\paragraph{Preorders of Hom-Sets. }In our case, we desire preorders \textit{within }hom-sets, $\displaystyle \mathbf{C}( a,b)_{\leqslant } \cup 1$. We associate each hom-set $\displaystyle \mathbf{C}( a,b)$ with a pre-ordered set $\displaystyle \mathbf{C}( a,b)_{\leqslant } \cup 1$ where morphisms $\displaystyle f,g,\dotsc \in \mathbf{C}( a,b)$ are objects $\displaystyle f,g,\dotsc \in \mathbf{C}( a,b)_{\leqslant } \cup 1$ with transforms to superior elements and with an identifying morphism $\displaystyle f_{f}^{1} :1\rightarrow f$ for each member.

In a functor string diagram, we typically express $\displaystyle g:a\rightarrow b$ with the symbol ``$\displaystyle g$''. However, we can also identify it by the identifying morphism $\displaystyle g_{g}^{1} :1\rightarrow g$ in the category $\displaystyle \mathbf{C}( a,b)_{\leqslant } \cup 1$. Then, using $\displaystyle g_{g}^{1} =f_{f}^{1} ;\leqslant _{g}^{f}$, we can substitute this identifier $\displaystyle g_{g}^{1}$ with $\displaystyle f_{f}^{1}$ and the noted transformation $\displaystyle \leqslant _{g}^{f}$.

\begin{figure}[htb]
    \caption{Using a series of re-expressions, we can intuitively show superior morphisms as equivalent to a transform of another.}
    \label{fig:FSDpreorders}
    \includegraphics[scale=\scalemed]{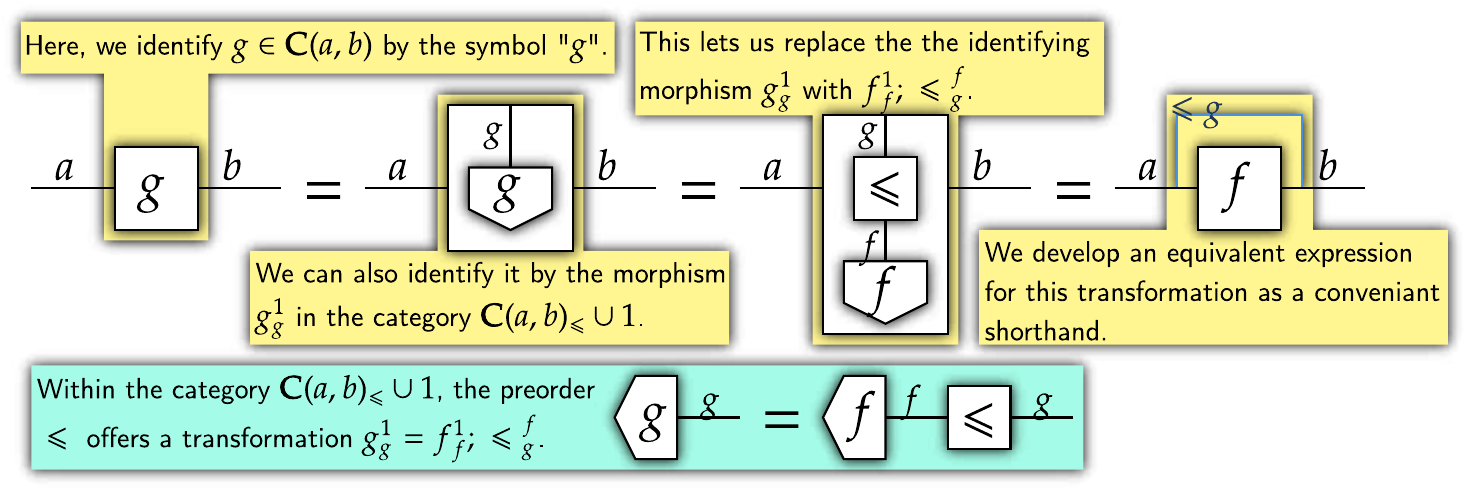}
\end{figure}

\paragraph{Two-Categories.} Consistent pre-orders on hom-sets in a category carry through. If $\displaystyle f\leqslant g$ for $\displaystyle f,g\in \mathbf{C}( a,b)$, and $\displaystyle k\leqslant h$ for $\displaystyle k,h\in \mathbf{C}( b,c)$, then we have $\displaystyle \leqslant _{g;h}^{f;k}$ in $\displaystyle \mathbf{C}( b,c)_{\leqslant }$. This requires $\displaystyle \mathbf{C}$ to have monotonic non-decreasing morphisms. Between adjacent pre-ordered hom-sets, we establish a bifunctor $\displaystyle \Uparrow _{abc} :\mathbf{C}( a,b)_{\leqslant } \cup 1\times \mathbf{C}( b,c)_{\leqslant } \cup 1\rightarrow \mathbf{C}( a,c)_{\leqslant } \cup 1$ where $\displaystyle \Uparrow ( 1,1) =1$, which implements composition of morphisms, and which is associative with respect to subsequent compositions.

\begin{figure}[htb]
    \floatbox[{\capbeside\thisfloatsetup{capbesideposition={left,top}}}]{figure}[\FBwidth]
    {\caption{We have composition between two-morphisms from adjacent hom-sets.}
    \label{fig:twocategory}}
    {\includegraphics[scale=\scalemed]{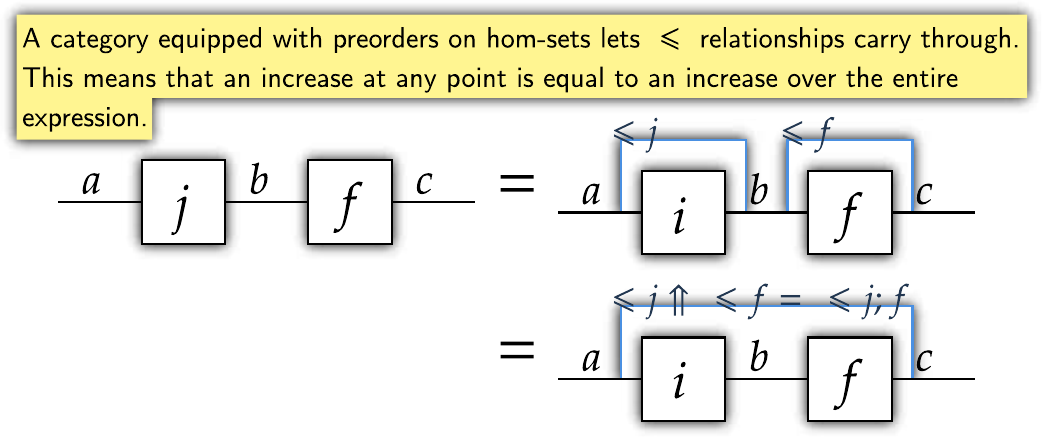}}
\end{figure}

Therefore, a consistent pre-order on hom-sets of a monotonic non-decreasing category implements the monotonic non-decreasing property that $\displaystyle f\leqslant g\rightarrow h;f\leqslant h;g$ categorically, by virtue of $\displaystyle \leqslant _{h}^{h} \Uparrow \leqslant _{g}^{f}$ yielding a transform $\displaystyle \leqslant _{h;g}^{h;f}$ in $\displaystyle \mathbf{C}( a,c)_{\leqslant } \cup 1$. This grants a two-category, and we show this in Figure \ref{fig:twocategory}.

Just as in algebraic manipulation of $\leqslant$, we are often more interested in the fact that a pre-order exists between expressions rather than the exact details of that preorder. So far, we have developed tools to view preorders as transforms, and have integrated these tools into the vertical sections of functor string diagrams. Now, we can confidently state lossy inequality expressions using diagrammatic categories. When we state $f \leqslant g$ for morphisms $a\rightarrow b$, we mean that there exists a transform $\leqslant^f_g$ between those morphisms in the preordered category $\textbf{C}(a,b)$.

\begin{figure}[htb]
    \caption{We use $\leqslant$ to state that \textit{a} transform exists. This lets clearly and readily manipulate various algebraic expressions.}
    \label{fig:preorder-short}
    \includegraphics[scale=\scalemed]{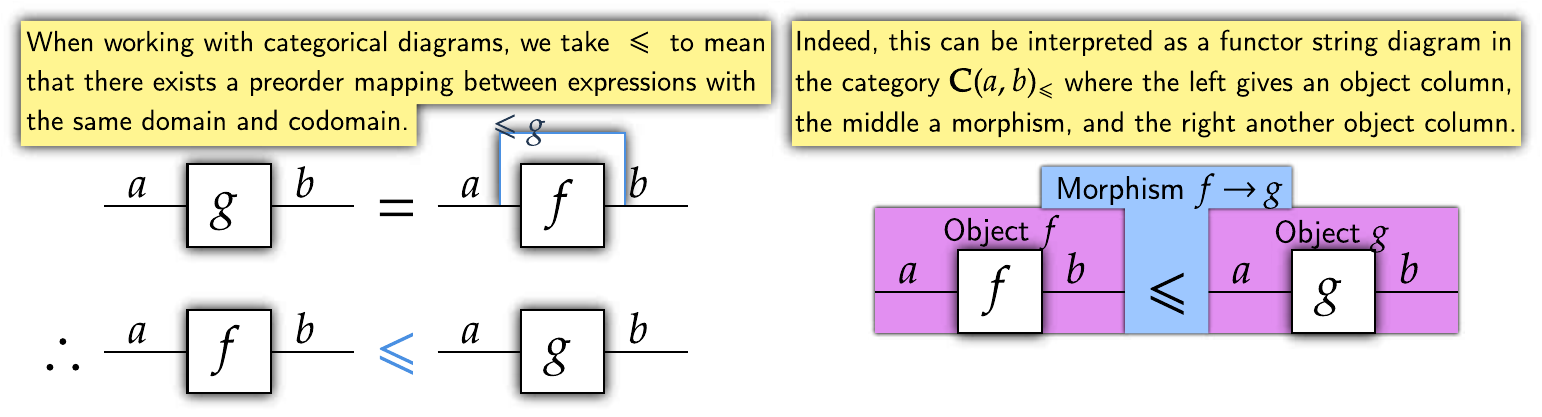}
\end{figure}

\subsection{Norphisms as Preorders}
\begin{figure}[htb] \floatbox[{\capbeside\thisfloatsetup{capbesideposition={left,top}}}]{figure}[\FBwidth]
    {\caption{Using functor string diagrams and preorders, we give categorical meaning to the equivariance condition as shown in Figure \ref{fig:inexact-composition}.}
    \label{fig:inexactpreorder}}
    {\includegraphics[scale=\scalemed]{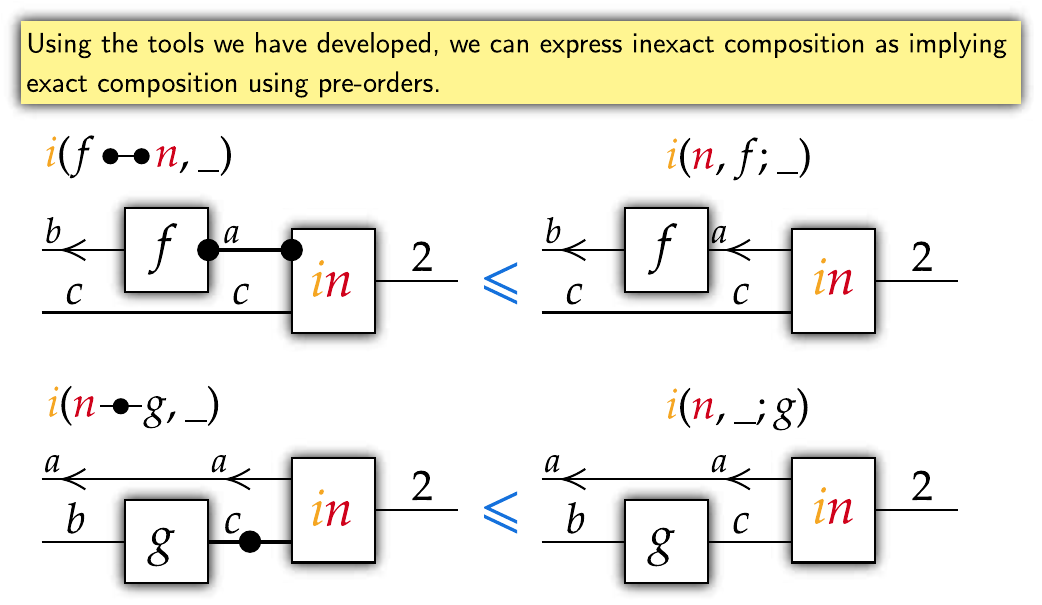}}
\end{figure}

We can now get to describing norphisms as purely categorical constructs by encoding the equivariance condition as a preorder as in Figure \ref{fig:inexactpreorder}. Furthermore, this preorder presentation encodes for the additional rule that bans described by norphisms are expansive, banning a set of morphisms and all others which are superior. This encodes for $\displaystyle f$ being infeasible implying that $\displaystyle g$ is as well, should $\displaystyle f\leqslant g$, limiting the search space of un-banned morphisms. Therefore, we not only present norphisms in a more streamlined traditional categorical manner, but also present them more formally.

First, recognize that the set $\displaystyle 2_{\leqslant }$ has a pre-order structure $\displaystyle 0\leqslant 1$, which is akin to implication. Thus, we can simply state the equivariance of inexact composition as a pre-order transform. Furthermore, as this is an expression in a pre-ordered category with monotonic non-decreasing morphisms, the expansiveness property of pre-orders is naturally implemented.

Therefore, we have captured all the structure of categorical negative information diagrammatically.

\section{Applications of Diagrammatic Norphisms}
In addition to streamlining the theory of negative information by employing intuitive diagrams, diagrammatic category theory enables theorems to be intuitively seen enhancing our understanding during applications. Here, we will cover cases from the original nategories paper, including negative information stemming from the triangle inequality and co-designs, and show how diagrams can provide greater insight while offering a platform for clear communication which motivates wider adoption.

\subsection{The Triangle Inequality}
Triangle inequalities abound in many cases, as knowledge of a set of processes to be completed almost necessarily indicates that some holistic approach may exist which performs them collectively. This yields relationships of the form $\displaystyle L( f;g) \leqslant L( f) +L( g)$ in many cases, such as the processing of probabilistic grammars \cite{klein-parsing-2003}, or the construction of neural networks, as two layers trained separately will yield a worse loss than training them together \cite{he-deep-2015, he-identity-2016}. Therefore, these compositional constructs are of particular interest.

Furthermore, we often have concepts of minimum distance. For example, in path-finding we have a Euclidean ``as the crow flies'' distance which is necessarily shorter than all possible paths. These underestimates or admissible heuristics are of particular interest to A* search or Dijkstra's algorithm and can yield a generalized means to solve problems in a variety of fields. This generality of composition indicates that category theory would be an appropriate method of analysis, however, these applications rely on understanding negative information and how the \textit{inability }to form a composition propagates through composition rules. Hence, typical categorical approaches are insufficient and, instead, nategories are necessary. We need a method of considering how a minimum on one hom-set propagates given knowledge of possible paths in others.

In this section, we show how diagrammatic category theory can clearly dissect the algebra of triangle inequalities and subsequently yield a form of norphism composition. The brevity of this section highlights the utility of the tools we have developed, they allow for the fundamental positive and negative compositional structure of systems to be clearly seen, a significant advantage over current methods and a substantial inroad towards considering this class of applied compositional problems in a universal, categorical manner.

\begin{figure}[htb] \floatbox[{\capbeside\thisfloatsetup{capbesideposition={left,top}}}]{figure}[\FBwidth]
    {\caption{The triangle inequality can be rearranged categorically, yielding a $\leq$ expression. This, of course, can be viewed as a transform within the category $\mathbf{Set}(\mathbf{C}(a,b),\mathbb{R})_\leqslant \cup 1$.}
    \label{fig:triangle}}
    {\includegraphics[scale=\scalemed]{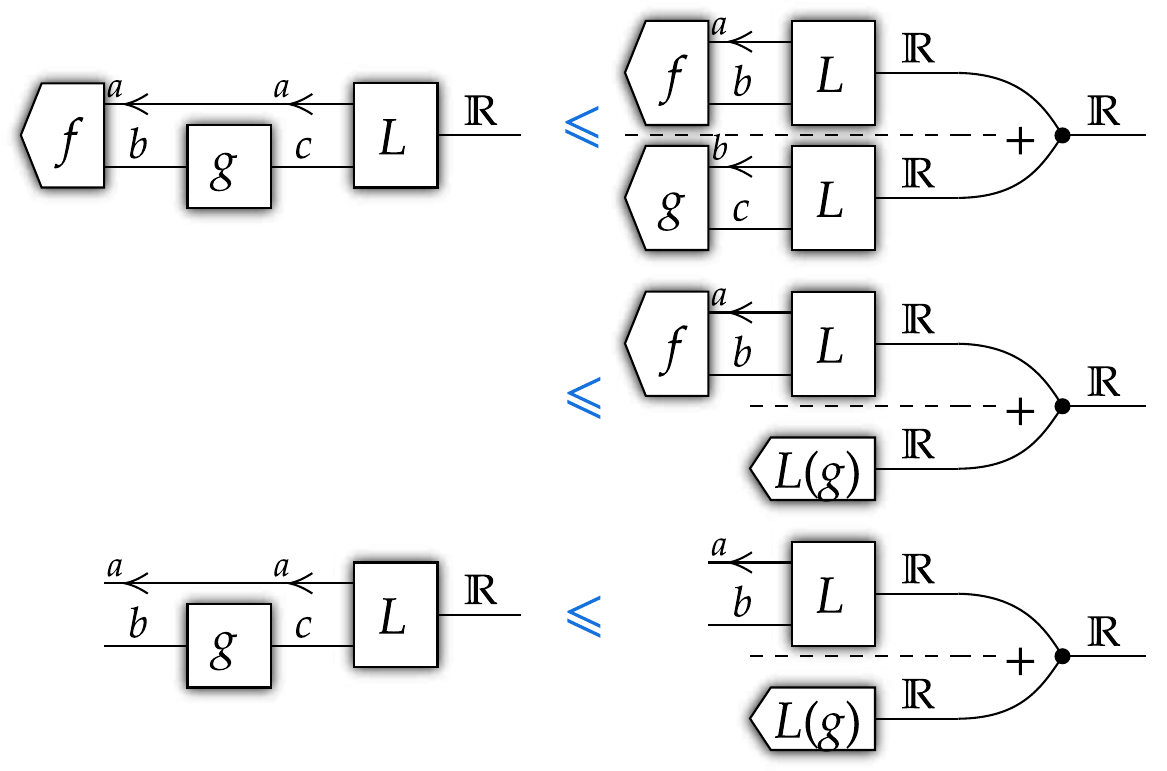}}
\end{figure}

We start by taking the triangle inequality --- $L(f;g) \leqslant L(f)+L(g)$ for some generally defined $L:\mathbf{C}(a,b)\rightarrow \mathbf{R}$ --- and considering it diagrammatically as in Figure \ref{fig:triangle}. This allows us to rearrange it so that we get an expression over all elements $f$ as morphisms of $\mathbf{Set}(1,\mathbf{C}(a,b))$. As this expression holds for all elements, it fully identifies the function in $\mathbf{Set}$ which operates on.

\begin{figure}[htb]
    \caption{We can define inexact composition between a morphism $g:b\rightarrow b$ and $\mathbf{C}(a,c)\rightarrow \mathbf{R}$ as we have a $\leqslant$ relationship.}
    \label{fig:inexact-R}
    \includegraphics[scale=\scalemed]{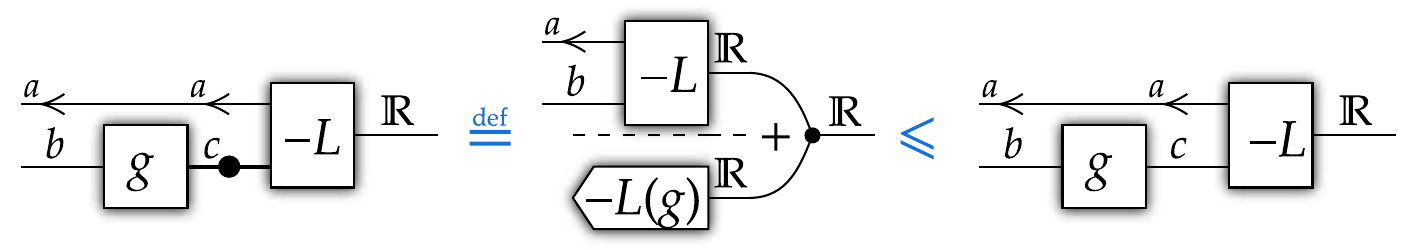}
\end{figure}

We can take the negative of both sides to arrive at something with a similar form to norphism inexact compositions shown in Figure \ref{fig:inexactpreorder}. In Figure \ref{fig:inexact-R}, we see that we can define inexact composition with a domain object $\displaystyle \mathbb{R}$ rather than $\displaystyle 2$. Nonetheless, this expression is clearly defined and is analogous to the categorical structure we have established so far. We see here an advantage of the general approach we have used, allowing negative information to be generalized.

\begin{figure}[!htb] \floatbox[{\capbeside\thisfloatsetup{capbesideposition={left,top}}}]{figure}[\FBwidth]
    {\caption{We can modify the inexact composition in Figure \ref{fig:inexact-R} to give inexact composition with codomain object $2$ by using the monotonic non-decreasing morphism $\leq -\mu$. This gives inexact composition mapping to $2$, yielding norphism composition.}
    \label{fig:inexact-2}}
    {\includegraphics[scale=\scalemed]{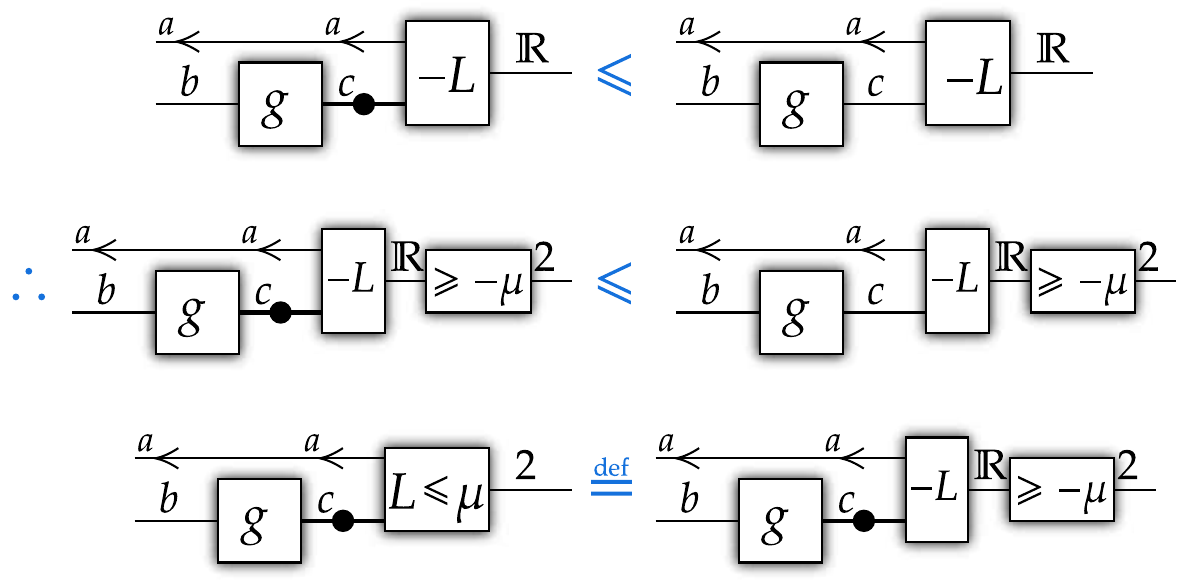}}
\end{figure}

Within a consistent pre-ordered category with monotonic non-decreasing morphisms, we can safely compose additional morphisms onto the above expressions, preserving the pre-order structure. A lower bound $\displaystyle \mu $ on the distance $\displaystyle L$ between objects $\displaystyle a$ and $\displaystyle c$ can be implemented by a norphism $\displaystyle ( -L) ;( \geqslant -\mu ) :\mathbf{C}( a,c)\rightarrow 2$. We can inherit the inexact composition from Figure \ref{fig:inexact-R} as shown in Figure \ref{fig:inexact-2}, deriving a norphism in a category with a triangle inequality.

With this algebra, we clearly see how the $\displaystyle \leqslant $ expression related to the triangle inequality in Figure~\ref{fig:triangle} becomes associated to a norphism by composition with monotonic non-decreasing morphisms in $\displaystyle \mathbf{Set}_{\leqslant }$ as in Figure~\ref{fig:inexact-2}. The tools we have used so far are all purely categorical and largely straightforward. Nonetheless, we are able to derive a norphism and the propagation of negative information concerning some ban on morphisms.

This shows the applied utility of the tools we have developed, and their potential to be further developed to be applied to a host of domains. Furthermore, we see that negative information is--indeed--purely categorical, emphasizing the usefulness of category theory even for considering what may seem to be atypical forms of composition.

\section{Diagrammatic Co-designs} \label{sec:diagrammatic-co-designs}
Complex systems co-design is an exciting avenue for applied category theory~\cite{censi2024, Zardini2023}. 
Interestingly, they have powerful diagrammatic semantics.
They describe the relationship between resources and processes, encoding how functionality demands relate to provided resources.
This can be applied to solve a wide range of design optimization problems, with applications in autonomy, mobility, and automotive~\cite{zardini2021ecc,zardini2021iros,ZardiniTask22,zardini2023camod,neumann-24-iv}.

Co-designs are a domain where considering negative information is critical. As they relate to the available resources and demanded functionalities of a system, we often want to utilize the knowledge that some process is known to be impossible. A basic example is the knowledge that physical mass or energy cannot be created. These limitations are norphisms --- a restriction on the morphisms $\displaystyle \mathbf{P}\rightarrow \mathbf{Q}$. These limitations have flow on affects to other objects as, otherwise, the bans they impose can be circumvented. This negative information can be considered with traditional nategories by establishing special composition rules.

However, we can show that co-designs accept negative information in a far more natural, and more categorical, manner. Indeed, by employing the insights of diagrams, we find that co-designs accept internal norphisms, meaning that all bans are represented as morphisms within the category itself. This insight allows us to use traditional codesign solving tools to investigate and consider negative information. %

\begin{figure}[htb] \floatbox[{\capbeside\thisfloatsetup{capbesideposition={left,top}}}]{figure}[\FBwidth]
    {\caption{Just as in monoidal string diagrams \cite{selinger-survey-2011}, we do not draw the unit object $\displaystyle \mathbf{1}$, corresponding to the Boolean vector space with a single value. Functionality requirements are morphisms $\displaystyle f:\mathbf{1}\rightarrow \mathbf{P}$ and resource availabilities are morphisms $\displaystyle r:\mathbf{P}\rightarrow \mathbf{1}$. Contraction indicates whether the requirements are satisfied by the available resources. This is also satisfied by superior functionality being satisfied by inferior resources.}
    \label{fig:vectorcovector}}
    {\includegraphics[scale=\scalemed]{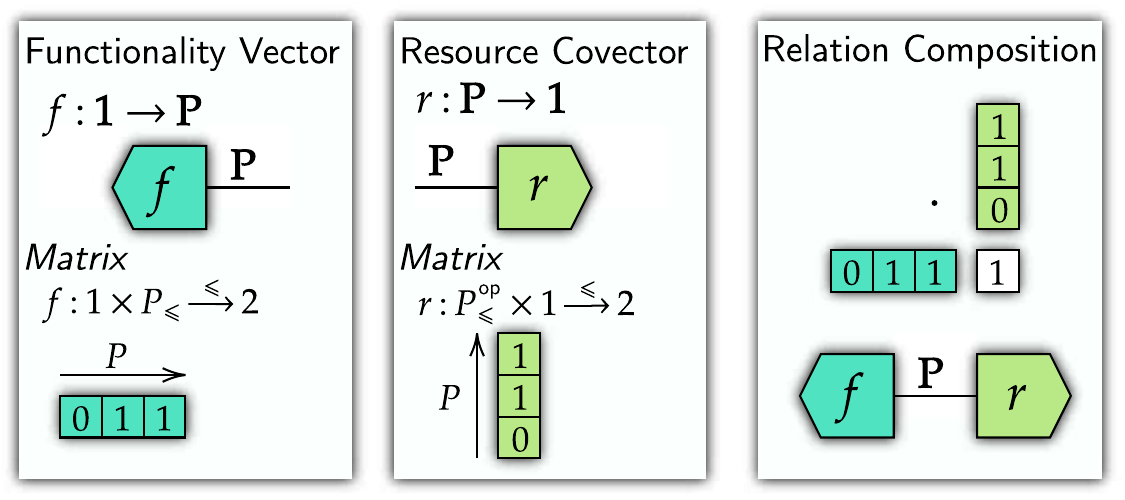}}
\end{figure}

Co-designs $\displaystyle \mathbf{DP}$ is a relations category where objects $\displaystyle \mathbf{P}$ are Boolean vector spaces indexed by the set $\displaystyle P$ and morphisms $\displaystyle d$ are relations between vector spaces, employing typical relation contraction as composition. $\displaystyle \mathbf{DP}$ is distinct, however, in that the Boolean matrices which correspond to relations have a particular pre-ordered structure as shown in Figure \ref{fig:vectorcovector}. Vectors, objects corresponding to output codomains, represent functionality requirements and have a non-decreasing structure. This represents some functionality requirement being satisfied by superior resources. Covectors, objects corresponding to input domains, represent resource availabilities and have a non-increasing structure. The contraction of vectors and covectors, therefore, answers whether some functionality requirement (or more) is found in the available pool of resources (or less).

\begin{figure}[htb] \floatbox[{\capbeside\thisfloatsetup{capbesideposition={left,top}}}]{figure}[\FBwidth]
    {\caption{Morphisms in $\textbf{DP}$ are co-designs, relations between Boolean vector spaces. They are ensured to provide valid functionality vectors and availability covectors by having a monotonic non-decreasing structure across the codomain axis and a monotonic non-increasing structure across the domain axis.}
    \label{fig:designproblem}}
    {\includegraphics[scale=\scalemed]{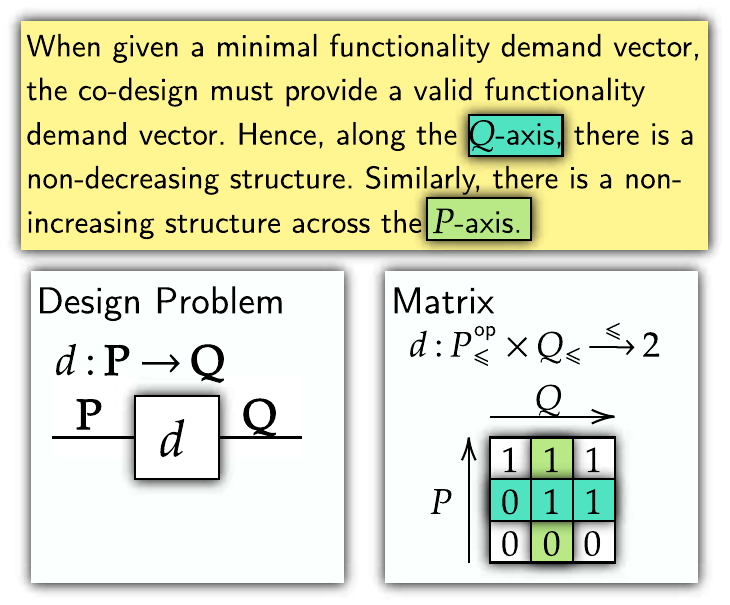}}
\end{figure}

Morphisms $\displaystyle d:\mathbf{P}\rightarrow \mathbf{Q}$ are design problems which process resources. Therefore, they map a functionality requirement in $\displaystyle P$, $\displaystyle f:\mathbf{1}\rightarrow \mathbf{P}$, to a functionality requirement in $\displaystyle \mathbf{Q}$, $\displaystyle f;d:\mathbf{1}\rightarrow \mathbf{Q}$. Similarly, they map a pool of available resources in $\displaystyle Q$, $\displaystyle r:\mathbf{Q}\rightarrow 1$, to a pool of available resources in $\displaystyle P$, $\displaystyle d;r:\mathbf{P}\rightarrow \mathbf{1}$. This is achieved by standard relation composition. These relations must, upon composition, produce valid non-decreasing vectors and valid non-increasing covectors. Therefore, their corresponding Boolean matrices must have a structure $\displaystyle d:P_{\leqslant }^{\text{op}} \times Q_{\leqslant }\xrightarrow{\leqslant } 2$ such that contraction along either axis yields the appropriate pre-ordered structure, as shown in Figure \ref{fig:designproblem}.

Relations can be understood as linear operations where addition is replaced by $\displaystyle \lor $ and multiplication by $\displaystyle \land $. We can understand the truth value of output indexes as whether a corresponding interaction between non-negative linear operations would yield non-negative values at those locations. From this, it follow that all the rules of manipulating linear operations are present. Furthermore, as both $\displaystyle \lor $ and $\displaystyle \land $ are monotonic non-decreasing, relations are always monotonic non-decreasing. Therefore, $\displaystyle \mathbf{DP}$, as a subcategory of relations, has all monotonic non-decreasing morphisms. As inputs (functionality demanded) increases, so does the output (available resources required).

\paragraph{Norphisms in Co-design.} 
The nature of norphisms in co-designs states impossibility results, e.g., that some functionality demand cannot be achieved by some pool of available resources. 
For instance, we can state that $\displaystyle f:\mathbf{1}\rightarrow \mathbf{Q}$ cannot be provided by $\displaystyle \varphi :\mathbf{P}\rightarrow \mathbf{1}$ (\textit{the objects }$\displaystyle \mathbf{1}$\textit{ correspond to Boolean vector spaces }$\displaystyle 2^{1}$) with any choice of design problem $\displaystyle d:\mathbf{P}\rightarrow \mathbf{Q}$.

\begin{figure}[htb] \floatbox[{\capbeside\thisfloatsetup{capbesideposition={left,top}}}]{figure}[\FBwidth]
    {\caption{For a locally small category $\displaystyle \mathbf{C}_{\leqslant }$ with a separating object, we can construct a performance norphism which bans all morphisms $\displaystyle x\rightarrow y$ which return $\displaystyle \varphi $ or greater for some input $\displaystyle \psi :1\rightarrow x$.}
    \label{fig:performancenorphism}}
    {\includegraphics[scale=\scalemed]{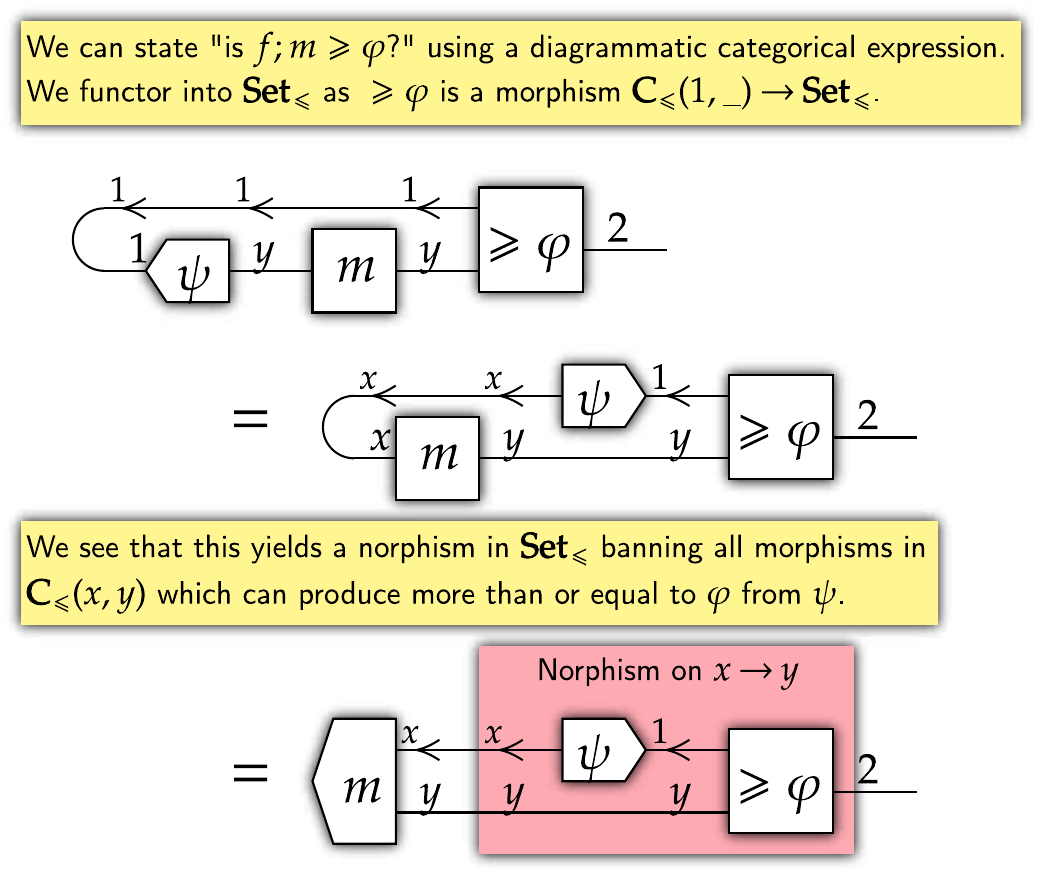}}
\end{figure}

Generally, for a locally small category $\displaystyle \mathbf{C}_{\leqslant }$ with a separating object $\displaystyle 1$ (meaning morphisms $\displaystyle x\rightarrow \_$ are uniquely identified by the composition with the set of morphisms $\displaystyle 1\rightarrow x$), we have a faithful hom-functor into set $\displaystyle \mathbf{C}_{\leqslant }( 1,\_) :\mathbf{C}_{\leqslant }\rightarrow \mathbf{Set}_{\leqslant }$. The preorder on $\displaystyle 1\rightarrow y$ can be used to generate a predicate $\displaystyle ( \geqslant \varphi ) :\mathbf{C}_{\leqslant }( 1,y)\rightarrow 2$ for some $\displaystyle \varphi :1\rightarrow y$. We can use this to test the statement ``is $\displaystyle \psi ;m\geqslant \varphi $?''. With the algebra of hom-functors \cite{abbott-robust-2023}, we rearrange this expression to get a generic ``performance norphism'' $\displaystyle \mathbf{C}_{\leqslant }( x,y)\rightarrow 2$ as in Figure \ref{fig:performancenorphism}.

Therefore, we see that norphisms of the type we desire which place some limit on performance can be naturally constructed from hom-functors and natural transformations. 
Co-designs are relations and, therefore, have linear structure as their objects are fundamentally vector spaces and their composition is fundamentally matrix multiplication with all positive values, albeit ignoring the magnitude of those values. Therefore, they are closed categories with access to internal hom-functors, transpositions, and monoidal products.

\begin{figure}[htb]
    \centering
    \includegraphics[scale=\scalemed]{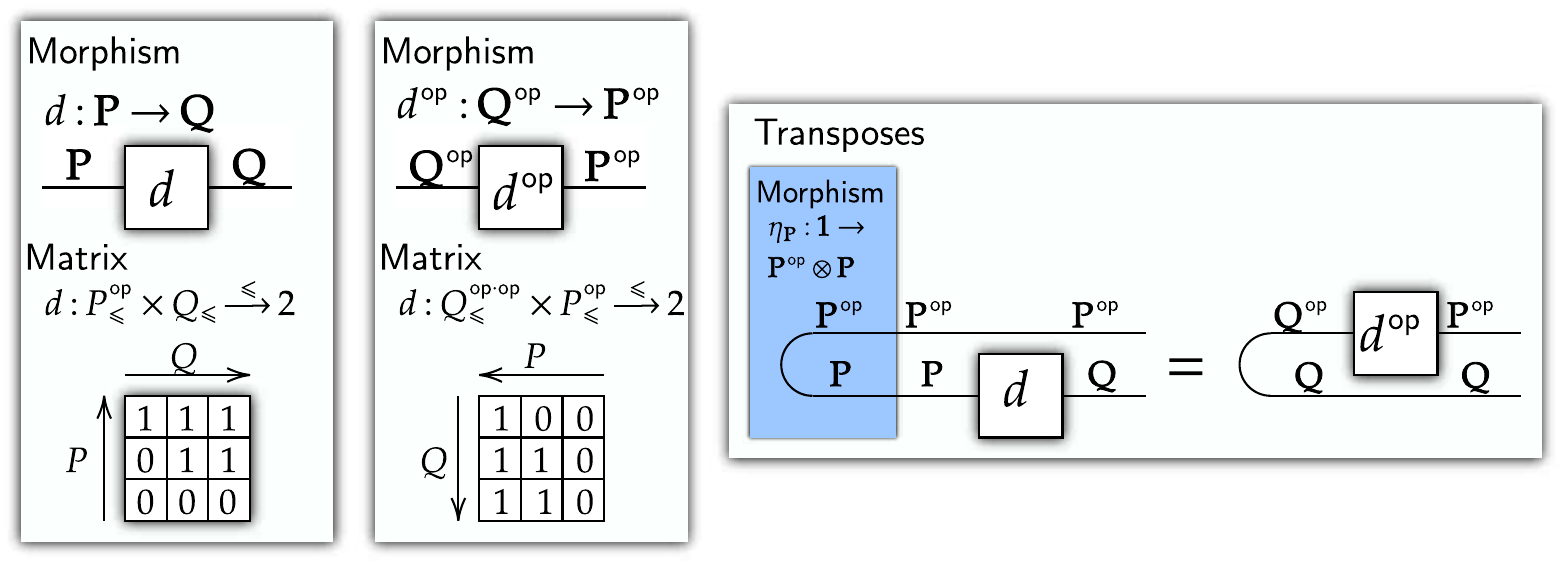}
    \caption{Within linear categories, we have internal hom-functors provided by outer $\displaystyle \lor $-products. Transposed design problems on product axes offer natural transformations.}
    \label{fig:transposes}
\end{figure}

However, we are required to fix types to ensure that the correct pre-order structure is present within the functionality demand vector spaces and resource availability provision covector spaces. This is achieved by having the monoidal unit $\displaystyle \eta _{\mathbf{P}} :1\rightarrow \mathbf{P}^{\text{op}} \otimes \mathbf{P}$, which provides transposes, switch $\displaystyle \mathbf{P}$ to a dual object $\displaystyle \mathbf{P}^{\text{op}}$ \cite{selinger-survey-2011}. By transposing and swapping the order of $\displaystyle P_\leqslant$, we get well-defined transposed design problems as shown in Figure \ref{fig:transposes}.

\begin{figure}[htb] \floatbox[{\capbeside\thisfloatsetup{capbesideposition={left,top}}}]{figure}[\FBwidth]
    {\caption{Relations offer all the tools to include a performance norphism internally.}
    \label{fig:codesignperformance}}
    {\includegraphics[scale=\scalemed]{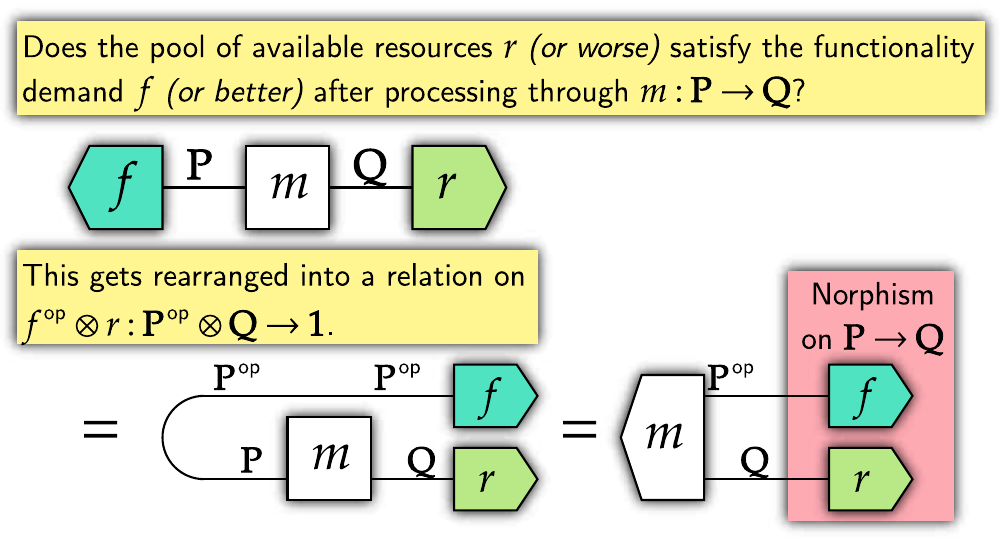}}
\end{figure}

Using (op)monoidal products for internal hom-functors and transposes for natural transformations, we can construct the performance norphism \textit{within }$\displaystyle \mathbf{DP}$ as show in Figure \ref{fig:codesignperformance}. The separating object in $\displaystyle \mathbf{DP}$ is $\displaystyle \mathbf{1}$, the vector space consisting of one Boolean value. The $\displaystyle \psi :1\rightarrow x$ morphism is $\displaystyle f:\mathbf{1}\rightarrow \mathbf{P}$ while the $\displaystyle \varphi :1\rightarrow y$ morphism against which we test is $\displaystyle r:\mathbf{Q}\rightarrow \mathbf{1}$. The performance norphism $\displaystyle \mathbf{C}_{\leqslant }( x,y)\rightarrow 2$ corresponds to a design problem $\displaystyle \mathbf{P}^{\text{op}} \otimes \mathbf{Q}\rightarrow \mathbf{1}$. Therefore, the performance norphism asks if $\displaystyle m:\mathbf{P}\rightarrow \mathbf{Q}$ (transposed into $\displaystyle \mathbf{1}\rightarrow \mathbf{P}^{\text{op}} \otimes \mathbf{Q}$) is able to achieve the functionality demand of $\displaystyle f$ from the available resources of $\displaystyle r$. Furthermore, as relations are monotonic non-decreasing, morphisms $\displaystyle m:\mathbf{P}\rightarrow \mathbf{Q}$ which provide more functionality than $\displaystyle f$ from $\displaystyle r$ are also banned.

\begin{figure}[htb] \floatbox[{\capbeside\thisfloatsetup{capbesideposition={left,top}}}]{figure}[\FBwidth]
    {\caption{Any design problem $\displaystyle \mathbf{P}^{\text{op}} \otimes \mathbf{Q}\rightarrow \mathbf{1}$ can be viewed as a norphism banning design problems $\displaystyle \mathbf{P}\rightarrow \mathbf{Q}$ which can provide functionality $\displaystyle f[i]$ from available resources $\displaystyle r[i]$ for any constituent $\displaystyle i$.}
    \label{fig:combocodesign}}
    {\includegraphics[scale=\scalemed]{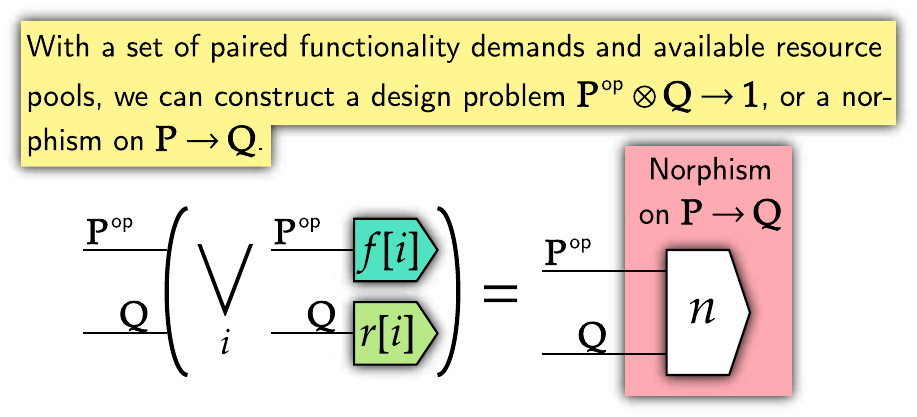}}
\end{figure}

Finally, as in Figure \ref{fig:combocodesign}, we can take a $\displaystyle \lor $-sum of different performance norphisms $\displaystyle n=\lor _{i}\left( f[ i]^{\text{op}} \otimes r[ i]\right)$ to get a relation $\displaystyle \mathbf{P}^{\text{op}} \otimes \mathbf{Q}\rightarrow \mathbf{1}$ which returns true if, for any $\displaystyle i$, the provided relation $\displaystyle m:\mathbf{P}\rightarrow \mathbf{Q}$ satisfies the functionality demand $\displaystyle f[ i]$ with the resources $\displaystyle r[ i]$. This process lets us construct any design problem $\displaystyle \mathbf{P}^{\text{op}} \otimes \mathbf{Q}\rightarrow \mathbf{1}$ as a combination of bans and to view these as norphisms.

\begin{figure}[htb] \floatbox[{\capbeside\thisfloatsetup{capbesideposition={left,top}}}]{figure}[\FBwidth]
    {\caption{Norphisms on design problems $\displaystyle P\rightarrow Q$ are themselves design problems $\displaystyle \mathbf{P}^{\text{op}} \otimes \mathbf{Q}\rightarrow \mathbf{1}$. They can be composed with design problems to form other design problems of the form $\displaystyle \_^{\text{op}} \otimes \_\rightarrow \mathbf{1}$, which correspond to other design problems, showing how negative information has a natural categorical structure.}
    \label{fig:codesigncomposition}}
    {\includegraphics[scale=\scalemed]{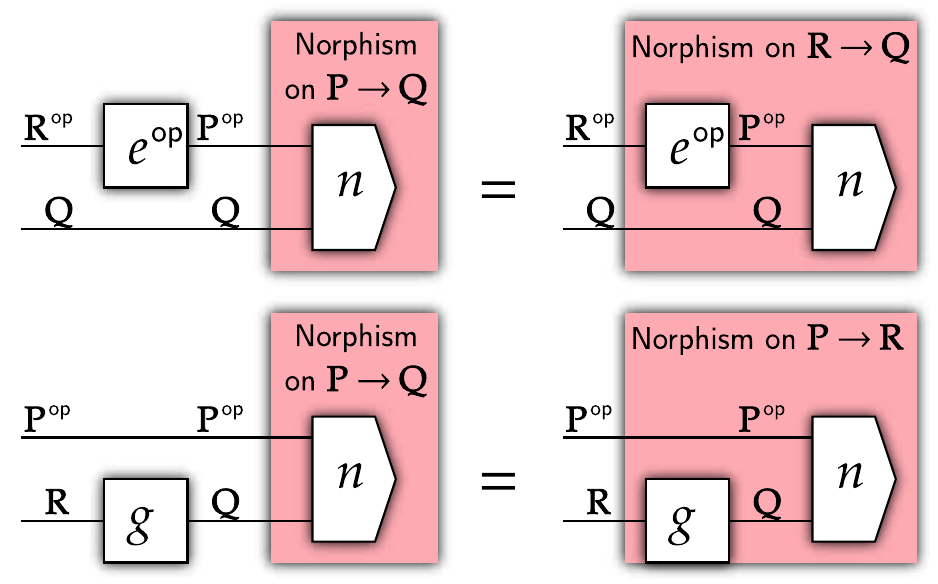}}
\end{figure}

The norphisms we have constructed exhibit propagation of negative information as this property follows from the axioms of categories. As displayed in Figure \ref{fig:codesigncomposition}, a norphism on $\displaystyle \mathbf{P}\rightarrow \mathbf{Q}$ is a design problem $\displaystyle n:\mathbf{P}^{\text{op}} \otimes \mathbf{Q}\rightarrow \mathbf{1}$ and provides a hook for composition with design problems $\displaystyle e:\mathbf{P}\rightarrow \mathbf{R}$ or $\displaystyle g:\mathbf{R}\rightarrow \mathbf{Q}$ to form norphisms on $\displaystyle \mathbf{R}\rightarrow \mathbf{Q}$ and $\displaystyle \mathbf{P}\rightarrow \mathbf{R}$ respectively. These represent the fact that a ban on $\displaystyle \mathbf{P}\rightarrow \mathbf{Q}$ and the existence of a morphism such as $\displaystyle e:\mathbf{P}\rightarrow \mathbf{R}$ implies a ban on $\displaystyle \mathbf{R}\rightarrow \mathbf{Q}$ must be placed, otherwise composition with $\displaystyle e$ would allow the ban to be circumvented. Norphisms composition in the case of co-designs, then, is exact.

\begin{figure}[htb] \floatbox[{\capbeside\thisfloatsetup{capbesideposition={left,top}}}]{figure}[\FBwidth]
    {\caption{We obtain a functionality vectors which corresponds to just more than what an availability covector provides by a transposed negation. As a performance norphism, this implements a physical resource being unable to be increased by any design problem.}
    \label{fig:schema}}
    {\includegraphics[scale=\scalemed]{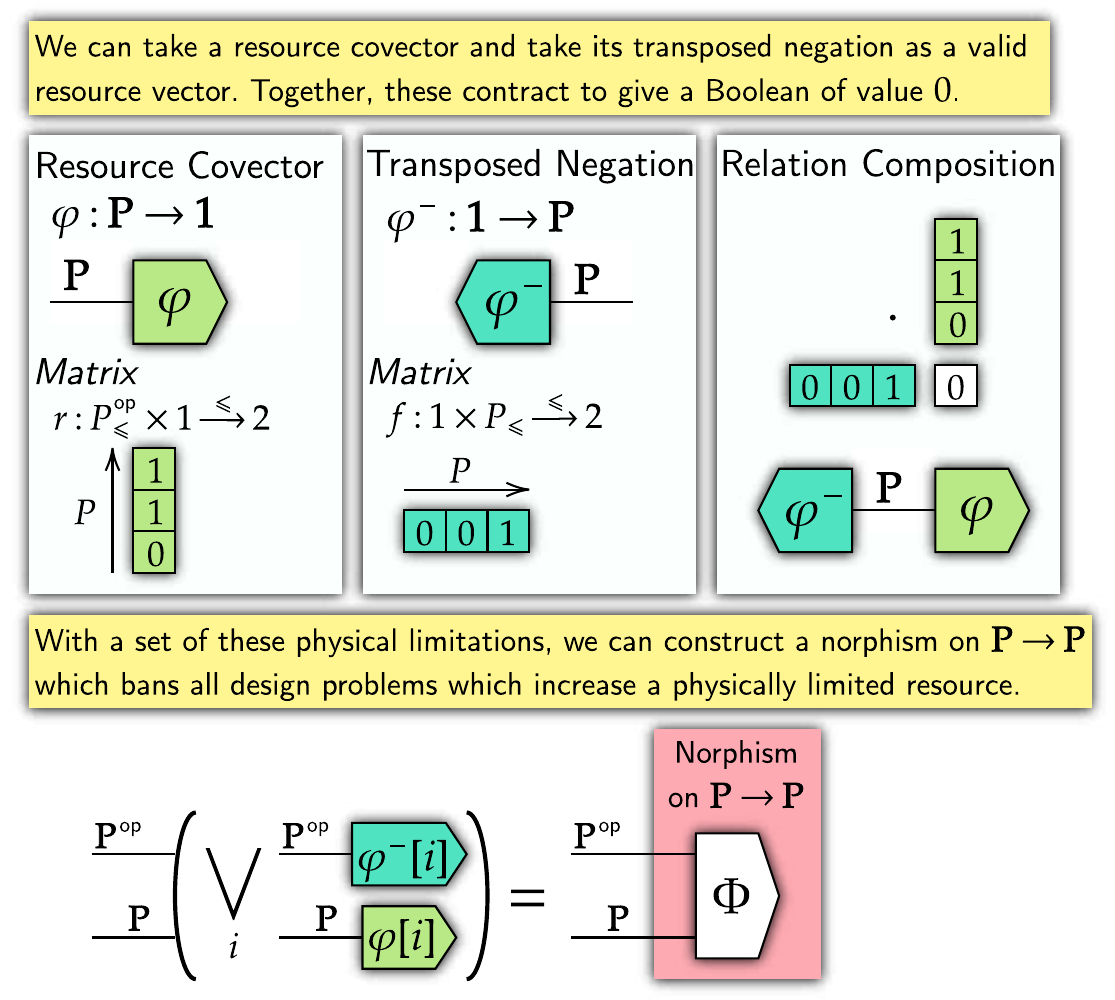}}
\end{figure}

\paragraph{Norphism Schemas. }This framework can also be used to implement physical resource constraints, wherein no composition of design problems can exist to expand a physical resource pool. Hence, a resource pool $\displaystyle \varphi :\mathbf{P}\rightarrow \mathbf{1}$ can never provide the functionality to recreate more than itself, $\displaystyle \varphi ^{+} :\mathbf{1}\rightarrow \mathbf{P}$. This can be directly achieved by the transposed negation converting the resource pool into a functionality pool, which we see in Figure \ref{fig:schema}. As there is no overlap between a pool and its transposed negation, this does not ban the identity. Instead, it bans all design problems which offer more than the identity between a resource and itself. We take a sum over all physically restricted components of the resource pool to implement a resource limitation norphism.

\section{Conclusion} Our diagrammatic investigation of co-designs and negative information reveals how basic categorical constructs explain a number of complex properties. Co-designs are shown to ultimately be relations with a fundamental linear structure, thereby inheriting graphical manipulation by monoidal string diagrams which provides intuitive insight into complex algebra. Furthermore, as we have previously shown that norphisms arise from predicates over hom-sets, the closedness of the codesign category allows negative information to be considered internally. This leads to further natural properties. In the context of co-designs, then, this investigation reveals that category theory is just as capable of considering positive as negative information.

\nocite{*}
\bibliographystyle{eptcs}
\bibliography{generic}

\end{document}

%% file: utils/preamble_act.tex
\usepackage{import}

\input{\rootdir/utils/packages_act}
\input{\rootdir/utils/tikz_utils}
\input{\rootdir/utils/symbols}
\input{\rootdir/utils/symbols_papers}
\input{\rootdir/utils/other_symbols}

%% file: utils/packages_act.tex
\usepackage[english]{babel} %

\usepackage[usenames,dvipsnames,table]{xcolor}

\usepackage{ifthen}
\usepackage{iftex}
\usepackage{ifxetex}
\usepackage[utf8]{inputenc}
\usepackage[style=russian]{csquotes}

\usepackage{blkarray, array}

\usepackage[T1]{fontenc}
\usepackage{lipsum}
\usepackage{xspace}
\usepackage{amsmath}
\usepackage{amsthm}
\usepackage{amssymb} %
\usepackage{mathtools}

\usepackage{multirow}
\usepackage{prftree}

\usepackage{comment}
\usepackage{accents}

\usepackage{ebproof}
\usepackage{units}

\usepackage{todonotes}

\usepackage{pgfplots}
\pgfplotsset{compat=1.16}
\usepackage[framemethod=TikZ]{mdframed}
\usepackage{graphicx}
\usepackage{calc}

\usepackage{adjustbox}
\usepackage[section]{placeins}
\usepackage{eso-pic}
\usepackage{capt-of}
\usepackage{wrapfig}
\usepackage[normalem]{ulem}
\usepackage{stmaryrd}
\usepackage{pifont}
\usepackage{relsize}

\usepackage{makeidx}

\usepackage{longtable}

\usepackage{multicol}

\usepackage{transparent}
\usepackage{datetime2}
\usepackage{anyfontsize}
\usepackage{ifmtarg}

\usepackage{underscore}

\graphicspath{
    {\rootdir/figures/},
    {\rootdir/pics/},
    {\rootdir/slides/},
    {\rootdir/videos/generated/},
    {\rootdir/figures/reits/},
    {\rootdir/figures/firstpaper/},
    {\rootdir/figures/uncertainty/},
    {\rootdir/papers/uncertainty/},
    {},
    {../}}

\newboolean{statuscolors}
\newboolean{instructors}
\newboolean{showslides} %
\newboolean{showproofs}
\newboolean{debugimages}
\newboolean{cachepdf}
\newboolean{devel} %

\newboolean{codeexercises}
\newboolean{githublink}

\newcommand{\instructors}[1]{\ifbool{instructors}{%
        #1
    }{}}

\newcommand{\showslides}[1]{\ifbool{showslides}{#1}{}}

\newcommand{\showproofs}[1]{\ifbool{showproofs}{#1}{}}
\newcommand{\devel}[1]{\ifbool{devel}{%
        \clearpage
        \pagecolor{yellow!7}
        #1
        \clearpage
        \pagecolor{white}
    }{}
}
\newcommand{\toexclude}[1]{\ifbool{devel}{%
        \clearpage
        \pagecolor{gray}
        #1
        \clearpage
        \pagecolor{white}
    }{}
}

\newcommand{\codeexercises}[1]{%
    \ifbool{codeexercises}{%
        \FloatBarrier
        \clearpage
        \pagecolor{blue!7}
        #1
        \clearpage
        \pagecolor{white}
    }{}%
}

\newlength{\dpiconwidth}
\newcommand{\linkgithub}{
    \ifbool{githublink}{%
        \usebox{\urlbox}%
    }{}
}
\newcommand{\draftnotice}{
    \ifbool{githublink}{%
        \footnotesize You are reading a draft compiled on \DTMnow%
    }{}%
}

\let\oldsubseteq\subseteq
\renewcommand{\subseteq}{\mathrel{\oldsubseteq}\linebreak[0]}
\let\oldsupseteq\supseteq
\renewcommand{\supseteq}{\mathrel{\oldsupseteq}\linebreak[0]}
\let\oldsubset\subset
\renewcommand{\subset}{\mathrel{\oldsubset}\linebreak[0]}
\let\oldsupset\supset
\renewcommand{\supset}{\mathrel{\oldsupset}\linebreak[0]}

\usepackage{filemod}

\usepackage{trimclip}

\usepackage{ulem}
\usepackage{silence}

%% file: utils/tikz_utils.tex
\usepackage{tikz}

\usetikzlibrary{positioning,
  intersections,
  hobby,
  patterns,
  calc,
  decorations.pathmorphing,
  decorations.markings,
  shadows,
  shapes,
  cd,
  decorations.markings,
  positioning,
  arrows.meta,
  shapes,
  calc,
  fit,
  quotes,
backgrounds}

\usepackage{tikzpagenodes}

\tikzcdset{arrow style=tikz}
\tikzstyle{arrowmorphismstyle} = [-{Triangle[width=3.5pt,length=3.5pt]},draw=morphisms]
\tikzstyle{arrownorphismstyle} = [-{Triangle[width=3.5pt,length=3.5pt]},dashed, draw=norphisms]
\tikzstyle{arrowmorphismstylemapsto} = [|-{Triangle[width=3.5pt,length=3.5pt]},draw=morphisms]
\tikzstyle{arrowfunctorstyle} = [-{Triangle[width=3.5pt,length=3.5pt]},draw=functors]
\tikzstyle{block} = [draw=morphisms, rectangle, minimum height=2em, minimum width=3em, thick]
\tikzstyle{blockgroup} = [draw=morphisms, rectangle, minimum height=2em, minimum width=3em, dotted, thick]
\tikzstyle{resblock} = [draw=red, rectangle, minimum height=2em, minimum width=3em, thick]
\tikzstyle{compblock} = [draw=blue, rectangle, minimum height=2em, minimum width=3em, thick]
\tikzstyle{effblock} = [draw=black, rectangle, minimum height=2em, minimum width=3em, thick]
\tikzstyle{tupblock} = [draw=morphisms, rectangle, minimum height=2em, minimum width=3em, thick]
\tikzstyle{styleobjects} = [draw=objects, thick]
\tikzstyle{morblock} = [draw=morphisms, rectangle, minimum height=2em, minimum width=3em, thick, rounded corners=2]
\tikzstyle{norblock} = [draw=norphisms, rectangle, minimum height=2em, minimum width=3em, thick, rounded corners=2]

\tikzstyle{blockk} = [draw, rectangle, minimum height=2.5em, minimum width=3.5em]
\tikzstyle{block1} = [draw, rectangle, minimum height=1.5em, minimum width=2.5em]
\tikzstyle{blockDyn} = [draw, rectangle, minimum height=2.5em, minimum width=3.5em, align=center, inner sep=10pt, thick, fill=white, copy shadow={draw=black,fill=black,opacity=1,shadow xshift=0.5ex,shadow yshift=-0.5ex}]
\tikzstyle{blockAlg} = [draw, rectangle, minimum height=1.5em, minimum width=2.5em, align=center, inner sep=10pt, thick]
\tikzstyle{sum} = [draw,circle]
\tikzstyle{input} = [coordinate]
\tikzstyle{output} = [coordinate]
\tikzstyle{pinstyle} = [pin edge={to-,thin,black}]
\tikzcdset{every label/.append style={font=\normalsize}}

\tikzset{fcname/.store in =\fcname, fcname={}}
\tikzset{funame/.store in =\funame, funame={}}
\tikzset{rcname/.store in =\rcname, rcname={}}
\tikzset{runame/.store in =\runame, runame={}}
\tikzset{whereres/.store in =\whereres, whereres=0.5}
\tikzset{wherefun/.store in =\wherefun, wherefun=0.5}
\tikzset{relres/.store in =\relres, relres={right}}
\tikzset{relfun/.store in =\relfun, relfun={left}}
\tikzset{posres/.store in =\posres, posres=0.95}
\tikzset{posfun/.store in =\posfun, posfun=0.95}
\tikzset{loos/.store in =\loos, loos=2}
\tikzset{feedback/.store in =\feedback, feedback=0}
\tikzset{
  DP/.style={%
    label/.style={
      font=\everymath\expandafter{\the\everymath\scriptstyle},
      inner sep=5pt,
      node distance=2pt and -2pt},
    semithick,
    node distance=1 and 1,
    rconn/.style={color=white,opacity=0.0,postaction={decorate}, shorten <=3.2pt, shorten >= 0.8,
    decoration={markings,
    mark= at position 0 with {
      \coordinate (a);
    },
    mark=at position .5 with
        {
        \ifthenelse{\equal{\feedback}{1}}{\def\angleOut{-90}\def\angleIn{-90}}{\def\angleOut{0}\def\angleIn{180}}
        \coordinate (b);
        \draw[dashed,\dpred,opacity=1.0] (a) to[out=\angleOut,in=\angleIn,looseness=\loos]
        node[pos=\posres,\relres=\whereres mm,\dpred,opacity=1,fill=white,inner sep=1pt,outer sep=1pt]{\rcname} (b);
      },
      mark= at position 1 with
        {
        \ifthenelse{\equal{\feedback}{1}}{\def\angleOut{0}\def\angleIn{0}}{\def\angleOut{180}\def\angleIn{0}}
        \ifthenelse{\equal{\feedback}{1}}{\def\symbol{\succeq}}{\def\symbol{\preceq}}
        \coordinate (c);
        \draw[\dpgreen,opacity=1.0] (c) to[out=\angleOut,in=\angleIn,looseness=\loos]
        node[pos=\posfun,\relfun=\wherefun mm,\dpgreen,opacity=1,fill=white,inner sep=1pt,outer sep=1pt]{\fcname} (b){}; %
        \node[draw,circle,inner sep=0.5pt,color=black,fill=white,opacity=1.0,line width=0.2mm] at (b) (nodepreceq) {$\symbol$};
      }
    }},
    runconn/.style={color=\dpred,dashed,postaction={decorate},
    decoration={markings,
    mark= at position 1 with {
      \coordinate (a);
      \draw[\dpred,opacity=1.0,dashed] ($(a)+(0.05,0)$) --++ (0.5,0) node[\relres,pos=\posres]{\runame};}
    }
    },
    funconn/.style={color=white,postaction={decorate},
    decoration={markings,
    mark= at position 0 with {
      \coordinate (a);
      \draw[\dpgreen] ($(a)+(-0.05,0)$) -- ($(a)+(-0.5,0)$) node[\relfun,pos=\posfun]{\funame};}
    }
    },
    execute at begin picture={\tikzset{
      x=\dpx, y=\dpy,
      every fit/.style={inner xsep=\dpx, inner ysep=\dpy}}}
  },
  dpx/.store in=\dpx,
  dpx = 1.5cm,
  dpy/.store in=\dpy,
  dpy = 1.5ex,
  dp port sep/.store in=\dpportsep,
  dp port sep=2,
  dp port length/.store in=\dpportlen,
  dp port length=4pt,
  dp min width/.store in=\dpminwidth,
  dp min width=0.5cm,
  dp rounded corners/.store in=\dpcorners,
  dp rounded corners=2pt,
  dp small/.style={dp port sep=1, dp port length=2.5pt, dpx=.4cm, dp min width=.4cm, dpy=.7ex},
  dp/.code 2 args={%
    \pgfmathsetlengthmacro{\dpheight}{\dpportsep * (max(#1,#2)) * \dpy}
    \pgfkeysalso{draw,%
      minimum width=\dpminwidth,%
      minimum height=\dpheight,%
      outer sep=0pt,%
      inner sep=5pt,%
      rounded corners=\dpcorners,
      thick,
      prefix after command={\pgfextra{\let\fixname\tikzlastnode}},
      append after command={\pgfextra{\draw
      \ifnum #1=0{} \else foreach \i in {1,...,#1} {
        ($(\fixname.north west)!{\i/(#1+1)}!(\fixname.south west)$) +(0,0) node[solid,left,circle,color=\dpgreen,draw,fill=\dpgreen,scale=0.3] {} coordinate (\fixname_fun\i) -- +(0,0) coordinate (\fixname_fun\i')}\fi %
      \ifnum #2=0{} \else foreach \i in {1,...,#2} {
        ($(\fixname.north east)!{\i/(#2+1)}!(\fixname.south east)$) +(0,0) coordinate (\fixname_res\i') -- +(0,0) node[solid,right,circle,color=\dpred,draw,fill=\dpred,scale=0.3] {} coordinate (\fixname_res\i)}\fi;
      }}}
  },
  dp name/.style={append after command={\pgfextra{\node[label=center,inner sep=2pt,fill=white] at (\fixname) {#1};}}}
}

\tikzset{
  oriented WD/.style={%
    every to/.style={out=0,in=180,draw},
    label/.style={
      font=\everymath\expandafter{\the\everymath\scriptstyle},
      inner sep=0pt,
      node distance=2pt and -2pt},
    semithick,
    node distance=1 and 1,
    decoration={markings, mark=at position .5 with {\arrow{stealth};}},
    ar/.style={postaction={decorate}},
    execute at begin picture={\tikzset{
      x=\bbx, y=\bby,
      every fit/.style={inner xsep=\bbx, inner ysep=\bby}}}
  },
  bbx/.store in=\bbx,
  bbx = 1.5cm,
  bby/.store in=\bby,
  bby = 1.75ex,
  bb port sep/.store in=\bbportsep,
  bb port sep=2,
  bb port length/.store in=\bbportlen,
  bb port length=4pt,
  bb min width/.store in=\bbminwidth,
  bb min width=1cm,
  bb rounded corners/.store in=\bbcorners,
  bb rounded corners=2pt,
  bb small/.style={bb port sep=1, bb port length=2.5pt, bbx=.4cm, bb min width=.4cm, bby=.7ex},
  bb/.code 2 args={%
    \pgfmathsetlengthmacro{\bbheight}{\bbportsep * (max(#1,#2)+1) * \bby}
    \pgfkeysalso{draw,minimum height=\bbheight,minimum width=\bbminwidth,outer sep=0pt,
      rounded corners=\bbcorners,thick,
      prefix after command={\pgfextra{\let\fixname\tikzlastnode}},
      append after command={\pgfextra{\draw
      \ifnum #1=0{} \else foreach \i in {1,...,#1} {
        ($(\fixname.north west)!{\i/(#1+1)}!(\fixname.south west)$) +(-\bbportlen,0) coordinate (\fixname_in\i) -- +(\bbportlen,0) coordinate (\fixname_in\i')}\fi %
      \ifnum #2=0{} \else foreach \i in {1,...,#2} {
        ($(\fixname.north east)!{\i/(#2+1)}!(\fixname.south east)$) +(-\bbportlen,0) coordinate (\fixname_out\i') -- +(\bbportlen,0) coordinate (\fixname_out\i)}\fi;
      }}}
  },
  bb name/.style={append after command={\pgfextra{\node[anchor=north] at (\fixname.north) {#1};}}}
}

\tikzset{
  tick/.style={postaction={
    decorate,
    decoration={markings, mark=at position 0.5 with {\draw[-] (0,.4ex) -- (0,-.4ex);}}}
  }
}

\newcommand{\saginclude}[1]{%
    \def\pdffile{\sagdir/\detokenize{#1}.pdf}%
    \def\srcfile{\sagdir/\detokenize{#1}.tikz}%
    \ifbool{cachepdf}{%
        \IfFileExists{\pdffile}{%
            \Filemodcmp{\pdffile}{\srcfile}{%
                \includegraphics[scale=1]{\pdffile}%
            }{%
                \PackageWarning{ACT4E}{TikZ file is more recent than outdated pdf cache \pdffile}%
                \input{\srcfile}%
            }%
        }{\input{\srcfile}}%
    }{%
        \input{\srcfile}%
    }%
}

\usetikzlibrary{calc}

\newlength{\brickwidth}
\newlength{\brickheight}
\newlength{\brickdia}
\newlength{\brickdiaheight}
\newlength{\brickmultipliedx}
\newlength{\brickmultipliedy}
\newlength{\halfbrickwidth}
\newlength{\hatchspread}
\newlength{\hatchthickness}
\newlength{\hatchshift}
\newcommand{\hatchcolor}{}
\tikzset{hatchspread/.code={\setlength{\hatchspread}{#1}},
  hatchthickness/.code={\setlength{\hatchthickness}{#1}},
  hatchshift/.code={\setlength{\hatchshift}{#1}},%
  hatchcolor/.code={\renewcommand{\hatchcolor}{#1}}}
\tikzset{hatchspread=3pt,
  hatchthickness=0.4pt,
  hatchshift=0pt,%
  hatchcolor=black}
\pgfdeclarepatternformonly[\hatchspread,\hatchthickness,\hatchshift,\hatchcolor]%
{custom north west lines}%
{\pgfqpoint{\dimexpr-2\hatchthickness}{\dimexpr-2\hatchthickness}}%
{\pgfqpoint{\dimexpr\hatchspread+2\hatchthickness}{\dimexpr\hatchspread+2\hatchthickness}}%
{\pgfqpoint{\dimexpr\hatchspread}{\dimexpr\hatchspread}}%
{%
  \pgfsetlinewidth{\hatchthickness}
  \pgfpathmoveto{\pgfqpoint{0pt}{\dimexpr\hatchspread+\hatchshift}}
  \pgfpathlineto{\pgfqpoint{\dimexpr\hatchspread+0.15pt+\hatchshift}{-0.15pt}}
  \ifdim \hatchshift > 0pt
  \pgfpathmoveto{\pgfqpoint{0pt}{\hatchshift}}
  \pgfpathlineto{\pgfqpoint{\dimexpr0.15pt+\hatchshift}{-0.15pt}}
  \fi
  \pgfsetstrokecolor{\hatchcolor}
  \pgfusepath{stroke}
}

\pgfdeclarelayer{background}
\pgfdeclarelayer{foreground}
\pgfsetlayers{background,main,foreground}

\newcommand{\rope}[7]{
\begin{tikzpicture}[n/.style={circle,draw,minimum size=2mm,inner sep=0, fill=#1}, remember picture]
\ifthenelse{\lengthtest{#5 pt=0 pt} \AND \lengthtest{#6 pt=1 pt}}
{\draw[ultra thick, #1] (0,0) -- node[pos=#5, #1,n]{} node[pos=#6,n] {} node[pos=0.5, above, black]{#7} (#3,#4){};}
{\draw[ultra thick, #1] (#5*#3,#5*#4) -- node[pos=0, #1,n,overlay]{} node[pos=1, #1,n] {} node[pos=0.5, above, black]{#7} (#6*#3,#6*#4){};
\begin{pgfonlayer}{background}
\draw[ultra thick, #2] (0,0)--(#5*#3,#5*#4){};
\draw[ultra thick, #2] (#6*#3,#6*#4)--(#3,#4){};
\end{pgfonlayer}}
\end{tikzpicture}}

\newcommand{\northenasymb}{%
\mathbin{\begin{tikzpicture}[n/.style={circle,draw,minimum size=0.8mm,inner sep=0, fill=black}]
\node[n] at (0.125,0.05) {};
\node at (0.1,0.1){};
\draw[thick](0,0.05)--(0.25,0.05){};
\end{tikzpicture}}}

\newcommand{\northenbsymb}{%
\mathbin{\begin{tikzpicture}[n/.style={circle,draw,minimum size=0.8mm,inner sep=0, fill=black}]
\node[n] at (0,0.05) {};
\node[n] at (0.2,0.05) {};
\node at (0.1,0.1){};
\draw[thick](0,0.05)--(0.2,0.05){};
\end{tikzpicture}}}

\newcommand{\northena}[2]{{#1 \northenasymb #2}}
\newcommand{\northenb}[2]{{#1 \northenbsymb  #2}}

\usepackage{forest}

\newcommand{\dirtree}[1]{%
    \begin{forest}
      for tree={
        font=\ttfamily,
        grow'=0,
        child anchor=west,
        parent anchor=south,
        anchor=west,
        calign=first,
        edge path={
          \noexpand\path [draw, \forestoption{edge}]
          (!u.south west) +(4.5pt,0) |- node[fill,inner sep=1.25pt] {} (.child anchor)\forestoption{edge label};
        },
        before typesetting nodes={
          if n=1
            {insert before={[,phantom]}}
            {}
        },
        fit=band,
        before computing xy={l=15pt},
      }
    #1
    \end{forest}
}

\definecolor{cblue}{rgb}{0.3,0.3,1.0}
\definecolor{cred}{rgb}{1.0,0.1,0.1}

\definecolor{cwhite}{rgb}{0.95,0.95,0.97}
\definecolor{funbox}{rgb}{0.92,1,0.9}
\definecolor{resbox}{rgb}{1,0.74,0.75}
\definecolor{impbox}{rgb}{1,0.92,0.85}

\usetikzlibrary{shadows.blur}
\usetikzlibrary{shapes.symbols}

\definecolor{red1}{rgb}{0.30,0,0}
\definecolor{red2}{rgb}{0.60,0,0}
\definecolor{red3}{rgb}{0.90,0,0}
\definecolor{my-cyan}{cmyk}{1,0,0,0}
\definecolor{my-magenta}{cmyk}{0,1,0,0}
\definecolor{my-yellow}{cmyk}{0,0,1,0}

\tikzcdset{
    setcd/.style={
        every matrix/.append style={
            draw=setcolorbord,
            fill=setcolor,
            rounded corners,
            inner sep=1pt,
            #1
        },
    },
}

\tikzcdset{
    posetcd/.style={
        every matrix/.append style={
            draw=posetcolorbord,
            fill=posetcolor,
            rounded corners,
            inner sep=1pt,
            #1
        },
    },
}

\tikzcdset{
    catcd/.style={
        every matrix/.append style={
            draw=catcolor,
            fill=catcolor,
            rounded corners,
            inner sep=4pt,
            #1
        },
    },
}

\tikzcdset{
    graphcd/.style={
        every matrix/.append style={
            draw=graphcolorbord,
            fill=graphcolor,
            rounded corners,
            inner sep=4pt,
            #1
        },
    },
}

%% file: utils/symbols.tex
\newcommand{\setinsymbol}{\in}
\newcommand{\setin}{\linebreak[0]\mathrel{{\color{formulasetcolor}\setinsymbol}}\linebreak[0]} %
\newcommand{\cartprod}{\mathbin{\stylesets{\raisedtexttimes}}} %
\definecolor{functionality}{rgb}{0.094869,0.500000,0.000000}
\definecolor{requirements}{rgb}{0.555789,0.000000,0.000000}
\definecolor{implementations}{RGB}{214,120,5}
\definecolor{functionalitylight}{rgb}{0.094869,0.800000,0.000000}
\definecolor{requirementslight}{rgb}{0.8,0.000000,0.000000}
\definecolor{implementationslight}{RGB}{214,120,5}
\definecolor{listcolor}{named}{black}
\definecolor{setcolor}{RGB}{243, 180, 126} %
\definecolor{setcolorbord}{RGB}{253, 143, 47}
\definecolor{posetcolor}{rgb}{0.8, 0.925, 0.992} %

\definecolor{posetcolorbord}{RGB}{81, 171, 252} %
\definecolor{graphcolor}{RGB}{252, 246, 207}%
\definecolor{graphcolorbord}{RGB}{219, 213, 173} %
\definecolor{formulasetcolor}{named}{brown} %
\definecolor{formulasetLcolor}{RGB}{213,87,96} %
\definecolor{elementscolor}{rgb}{0.1,0,0.4}
\definecolor{catcolor}{rgb}{0.94,1.0,0.94}
\definecolor{darkgreen}{rgb}{0.0, 0.5, 0.0}
\definecolor{darkred}{rgb}{0.5, 0.0, 0.0}
\definecolor{brightpink}{rgb}{1.0, 0.65, 0.79}
\definecolor{applegreen}{rgb}{0.55, 0.71, 0.0}
\definecolor{custompurple}{rgb}{0.18,0.46,0.74}
\definecolor{custompink}{rgb}{0.73,0.83,0.91}
\definecolor{funcname}{named}{blue}
\definecolor{classname}{named}{blue}
\definecolor{fieldname}{named}{darkgreen}
\definecolor{parula1}{RGB}{53,42,135}
\definecolor{parula2}{RGB}{15,92,221}
\definecolor{parula3}{RGB}{18,125,216}
\definecolor{parula4}{RGB}{7,156,207}
\definecolor{parula5}{RGB}{21,177,180}
\definecolor{parula6}{RGB}{21,177,180}
\definecolor{parula7}{RGB}{89,189,140}
\definecolor{parula8}{RGB}{165,190,107}
\definecolor{parula9}{RGB}{225,185,82}

\definecolor{dpcolor}{rgb}{0.5,0,0.5}
\def\dpgreen{functionality}
\def\dpred{requirements}

\definecolor{darkblue}{rgb}{0.0, 0.0, 0.55}
\definecolor{shadecolor}{rgb}{0.68, 0.85, 0.9}
\definecolor{transmuter}{RGB}{46,117,189}
\definecolor{objects}{RGB}{213,87,96}
\definecolor{morphisms}{RGB}{46,117,189}
\definecolor{norphisms}{RGB}{213,87,96}
\definecolor{colorposetsymbol}{RGB}{90,90,90}

\definecolor{transmuted}{RGB}{213,87,96}
\definecolor{naturaltransformations}{rgb}{0.1,0.5,0.0}
\definecolor{functors}{rgb}{0.5,0,0.5}

\definecolor{instructors}{RGB}{214,120,5}
\definecolor{upcolor}{RGB}{252, 3, 107}
\definecolor{staincola}{RGB}{181,23,0}
\definecolor{staincolb}{RGB}{255,150,141}
\definecolor{staincolc}{RGB}{181,23,30}
\definecolor{staincold}{RGB}{241,219,200}
\definecolor{staincole}{RGB}{92,24,13}
\definecolor{staincolf}{RGB}{29,177,0}
\definecolor{ropecola}{RGB}{187,127,14}
\definecolor{ropecolb}{RGB}{97,216,54}
\newcommand{\colTransmuted}{\color{transmuted}}
\newcommand{\styleobj}[1]{{\colTransmuted{#1}}}
\newcommand{\stylemorph}[1]{{\mathcolor{morphisms}{#1}}}
\newcommand{\stylenorph}[1]{{\mathcolor{norphisms}{#1}}}
\newcommand{\mthensymbol}{\fatsemi}
\newcommand{\mthen}{\mathbin{\stylemorph{\mthensymbol}}\linebreak[0]} %

\newcommand{\Obname}{Ob}
\newcommand{\Ob}{\styleobj{\mathrm{\Obname}}}

\newcommand{\CatSymbol}[1]{\textbf{#1}\xspace}
\newcommand{\Homname}{\textrm{Hom}}
\newcommand{\Hom}{\stylefunctors{\Homname}}

\newcommand{\HomSet}[3]{\Hom_{#1}\pars{{#2};\linebreak[0]{#3}}}
\newcommand{\CatC}{\CatSymbol{C}}

\newcommand{\Obja}{\styleobj{a}}

\newcommand{\Objb}{\styleobj{b}}

\newcommand{\Objc}{\styleobj{c}}

\newcommand{\nora}{\stylenorph{n}}

\newcommand{\mora}{\stylemorph{f}}
\newcommand{\morb}{\stylemorph{g}}

\def\bitcoinA{%
    \leavevmode
    \vtop{\offinterlineskip %
    \setbox0=\hbox{B}%
    \setbox2=\hbox to\wd0{\hfil\hskip-.03em
    \vrule height .3ex width .15ex\hskip .08em
    \vrule height .3ex width .15ex\hfil}
        \vbox{\copy2\box0}\box2}}

\definecolor{alitalia}{named}{darkgreen}
\definecolor{swiss}{named}{darkblue}
\definecolor{ewings}{named}{darkred}

\newlength{\ml}
\newlength{\mb}
\definecolor{bcolor}{rgb}{0.6,0.8,0.3}

\newcommand{\fatstart}{\llangle}
\newcommand{\fatend}{\rrangle}

\newcommand{\morab}{\mora\mthen\morb}

\newcommand{\cprod}{\mathbin{\stylemorph{\fatsemi}_{\hspace{-1pt}\stylemorph{\fatstart}}}} %

\newcommand{\tupconcatsymbol}{\fatsemi}
\newcommand{\tupconcat}{\mathbin{\stylemorph{\tupconcatsymbol}_{\hspace{-0.4pt}\stylemorph{\langle}}}} %

\newcommand{\Nom}{\stylefunctors{\mathrm{Nom}}}

\newcommand{\NomSet}[3]{\Nom_{#1}\pars{{#2};\linebreak[0]{#3}}}

%% file: utils/other_symbols.tex
\newcommand{\prf}{\expandafter\prftree[r]{\raisebox{3pt}{}}}
\newcommand{\prfdouble}{\expandafter\prftree[double]}
\newcommand{\prfcomma}{\expandafter\prftree[r]{\raisebox{3pt}{,}}}
\newcommand{\prfperiod}{\expandafter\prftree[r]{\raisebox{3pt}{.}}}
\newcommand{\prfsemi}{\expandafter\prftree[r]{\raisebox{3pt}{;}}}
\newcommand{\prfdoublecomma}{\expandafter\prftree[double][r]{\raisebox{3pt}{,}}}
\newcommand{\prfdoubleperiod}{\expandafter\prftree[double][r]{\raisebox{3pt}{.}}}

\makeatletter
\def\mathcolor#1#{\@mathcolor{#1}}
\def\@mathcolor#1#2#3{%
  \protect\leavevmode
  \begingroup
    \color#1{#2}#3%
  \endgroup
}
\makeatother

\DeclareMathSymbol{\shortminus}{\mathbin}{AMSa}{"39}
\newcommand{\minusone}{{{\shortminus}}\kern1pt1}

\makeatletter
\DeclareRobustCommand{\shortto}{%
    \mathrel{\mathpalette\short@to\relax}%
}

\newcommand{\short@to}[2]{%
    \mkern2mu
    \clipbox{{.4\width} 0 0 0}{$\m@th#1\vphantom{+}{\shortrightarrow}$}%
}
\makeatother

\makeatletter

\newcommand{\opensetbracket}{{\ColorIfNowBlack{formulasetcolor}{\mathbf{\{}}}}
\newcommand{\closesetbracket}{{\ColorIfNowBlack{formulasetcolor}{\mathbf{\}}}}}

\newcommand{\cleft}[2][.]{%
    \begingroup\colorlet{savedleftcolor}{.}%
    \color{#1}\mathopen{}\left#2\color{savedleftcolor}%
}
\newcommand{\cright}[2][.]{%
    \color{#1}\right#2\mathclose{}\endgroup%
}

\DeclareRobustCommand{\@makelists}[4]{%
    #1%
    \foreach \entry [count=\xi]  in {#4}{
        \ifnum\xi = 1\else#2\linebreak[0]\fi%
        {{\entry}}%
    }%
    #3%
}
\newcommand{\smallspace}{\mspace{4mu}}
\newcommand{\commaspace}{{,\smallspace}}
\DeclareRobustCommand{\makeset}[1]{%
    \@makelists{\opensetbracket}{\commaspace}{\closesetbracket}{#1}%
}
\DeclareRobustCommand{\makesetnocolor}[1]{%
    \@makelists{\{}{\commaspace}{\}}{#1}%
}
\DeclareRobustCommand{\tupp}[1]{%
    \@makelists{\langle}{\commaspace}{\rangle}{#1}%
}

\DeclareRobustCommand{\tup}[1]{%
        \@makelists{\mathopen{}\left\langle}{\commaspace}{\right\rangle\mathclose{}}{#1}%
}

\DeclareRobustCommand{\makecprod}[1]{%
    \@makelists{}{\cprod}{}{#1}%
}

\DeclareRobustCommand{\makecartprod}[1]{%
    \@makelists{}{\cartprod}{}{#1}%
}

\DeclareRobustCommand{\maketupconcat}[1]{%
    \@makelists{}{\tupconcat}{}{#1}%
}

\DeclareRobustCommand{\makesetBig}[1]{%
    \@makelists{{\mathopen{}\color{formulasetcolor}\mathbf{\Big\{}}}{\commaspace}{{\color{formulasetcolor}\mathbf{\Big\}}\mathclose{}}}{#1}%
}

\DeclareRobustCommand{\makesett}[1]{%
    \@makelists{\cleft[formulasetcolor]{\{}}{\commaspace}{\cright[formulasetcolor]{\}}}{#1}%
}

\DeclareRobustCommand{\makelist}[1]{%
    \@makelists{[}{\commaspace}{]}{#1}%
}

\DeclareRobustCommand{\cObj}[1]{%
    \@makelists{\fatstart}{\commaspace}{\fatend}{#1}%
}

\makeatother

\DeclarePairedDelimiter\mypars{\lparen}{\rparen}
\newcommand{\pars}[1]{\mypars*{#1}}

\edef\BlackColor{\csname\string\color@.\endcsname}

\makeatletter
\newcommand{\@CurrentColor}{}%
\DeclareRobustCommand{\ColorIfNowBlack}[2]{%
    \edef\@CurrentColor{\csname\string\color@.\endcsname}%
    \IfStrEq{\@CurrentColor}{\BlackColor}{%
        {\color{#1}#2}%
    }{%
        #2%
    }%
}
\makeatother

%% file: symbols.tex
\definecolor{norcolor}{RGB}{213,87,96}

\definecolor{GCQuestionsColor}{RGB}{213,87,96} %
\definecolor{GCAnswersColor}{RGB}{46,117,189} %
\definecolor{GCExpl}{named}{orange}
\definecolor{GCBool}{named}{black}
\definecolor{transmuted}{named}{black}

\renewcommand{\Hom}{\mathcolor{morphisms}{\Homname}}
\renewcommand{\Nom}{\mathcolor{norcolor}{\mathrm{Nom}}}

\newcommand{\nmincompat}[2]{\mathcolor{GCExpl}{i}_{#1#2}}

\definecolor{formulasetcolor}{named}{black}

\renewcommand{\mthen}{\mathbin{{\mthensymbol}}\linebreak[0]} %

\renewcommand{\cartprod}{\mathbin{\times}}